\documentclass{birkjour}
\usepackage{amssymb}
 \newtheorem{theorem}{Theorem}
 \newtheorem*{theorem*}{Theorem}
 \newtheorem{corollary}{Corollary}
 \newtheorem{lemma}{Lemma}
 \newtheorem*{lemma*}{Lemma}
 \newtheorem{statement}{Statement}
 \newtheorem*{statement*}{Statement}
 \theoremstyle{definition}
 \newtheorem{definition}{Definition}
 \theoremstyle{remark}
 \newtheorem{remark}{Remark}
 \newtheorem*{remark*}{Remark}
 
 \numberwithin{equation}{section}

\begin{document}

	\title[Miquel-Steiner's point locus]
		{Miquel-Steiner's point locus}

	\author[Yuriy Zakharyan]{Yuriy Zakharyan}

	\address{%
		Leninskie Gory 1\\
		119234 Moscow\\
		Russian Federation\\
		ORCID: 0000-0003-0042-3372
	}

	\email{yuri.zakharyan@gmail.com}

	\subjclass{51M05}

	\keywords{Miquel, Steiner, Brocard, Triangle, circle, cevians}

	\date{June 1, 2018}

	\begin{abstract}
		In this paper we reformulate Miquel-Steiner's theorem and we obtain Miquel-Steiner's point locus for an arbitrary triangle. We prove that this locus is related to conjugate circles and Brocard's circle. In addition, we obtain Miquel-Steiner's point locus in a case when cevians are perpendicular to each other, in a case when cevians form similar triangles. In addition, we prove that if Miquel-Steiner's point belongs to line, then cevinans intersection point belongs to line which is parallel to isogonal. Finally, we obtain few result for cases when Miquel-Steiner's point coincides with triangle centres.
	\end{abstract}

	\maketitle
	\section{Introduction}
		We start with a complete quadrilateral $AEFDBC$ ~[1, p.21] and with the Miquel-Steiner's point theorem ~[1, p.22].
		\begin{theorem*}[\textbf{Miquel, Steiner}]
			\label{th:miquel1}
			The circumcircles of four triangles of a complete quadrilateral are concurrent. Intersection point is called the Miquel-Steiner's point (Fig.1).
		\end{theorem*}
		\begin{center}
			\includegraphics[width=0.5\textwidth]{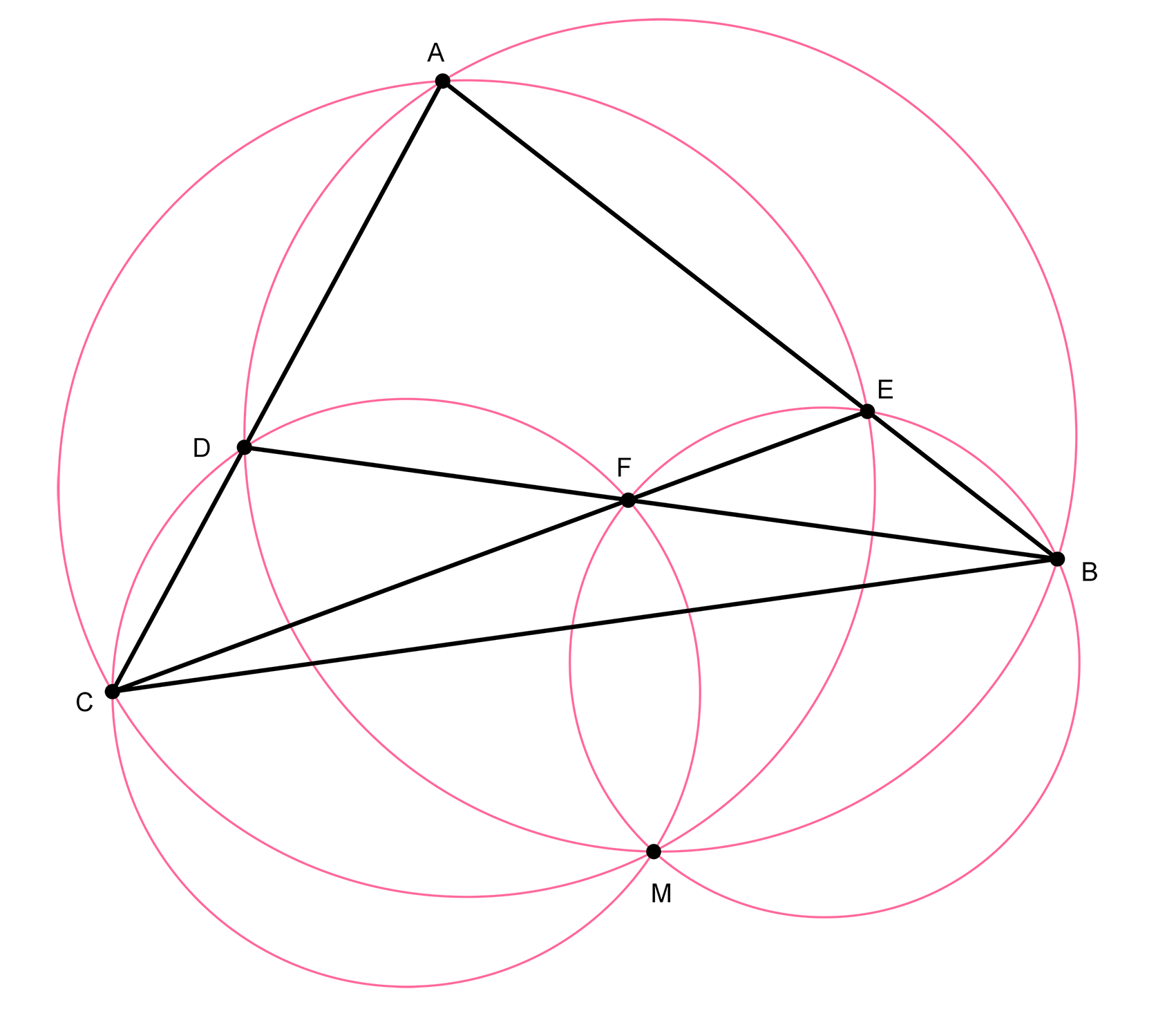}
			\newline
			Fig.1~Miquel-Steiner's point of complete quadrilateral
		\end{center}
		There is one narrowly known lemma, related to this theorem.
		\begin{lemma*}
			The Miquel-Steiner's point belongs to a diagonal of complete quadrilateral if and only if remaining quadrilateral vertices are cyclic (Fig.2).
		\end{lemma*}
		\begin{center}
			\includegraphics[width=0.5\textwidth]{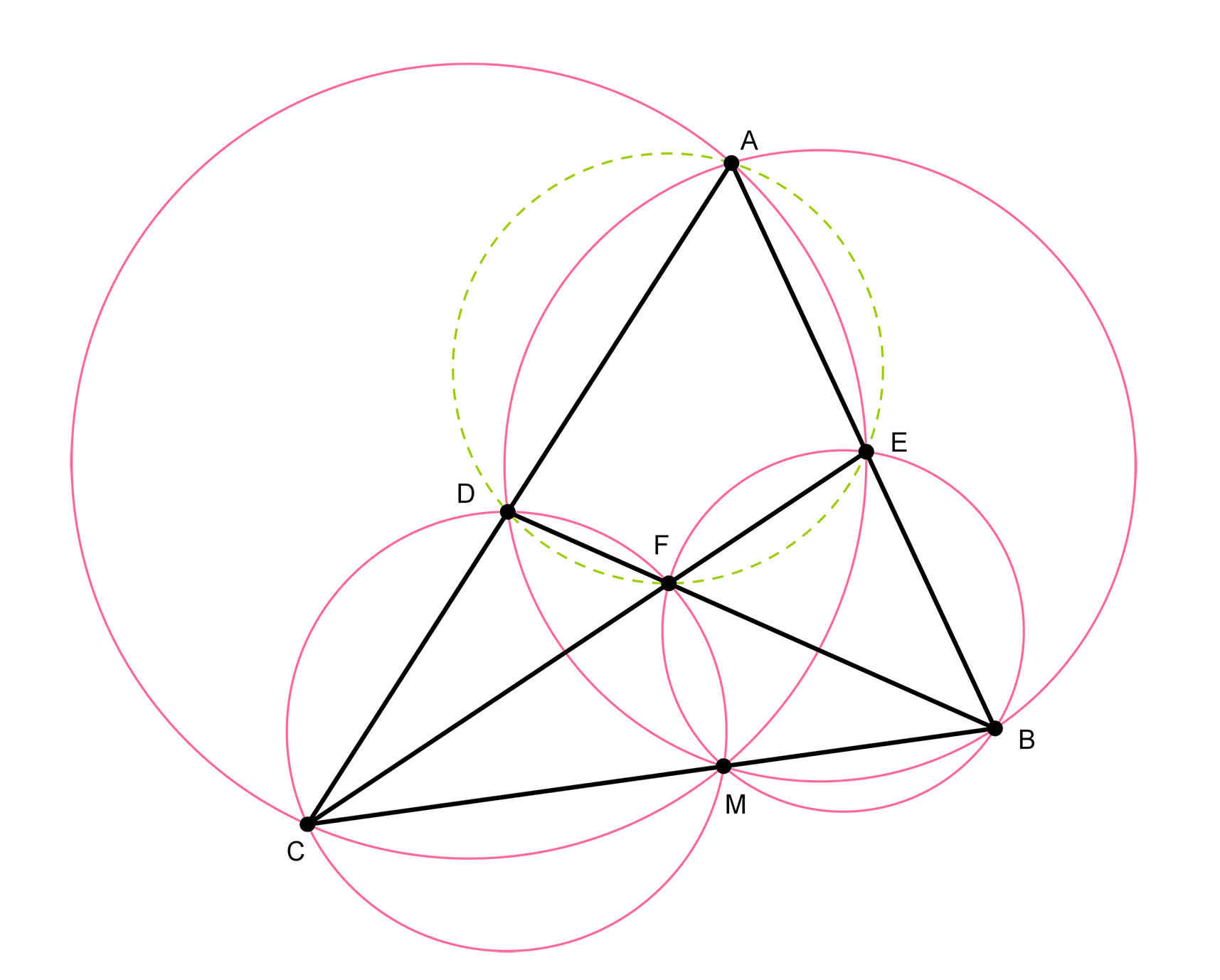}
			\newline
			Fig.2~Miquel-Steiner's point belongs to a diagonal
		\end{center}
		From now we start with arbitrary triangle $\Delta{ABC}$. Let us fix vertice $A$. Also let us take two cevians $BB_A, CC_A$ and their intersection point $N_A$. We can reformulate the Miquel-Steiner's point theorem.
		\begin{theorem*}
			The circumcircles of $\Delta{ABB_A}, \Delta{ACC_A}, \Delta{CB_AN_A}, \Delta{BC_AN_A}$  are concurrent. Intersection point $M_A$ is called the Miquel-Steiner's point (Fig.3).
		\end{theorem*}
		\begin{center}
			\includegraphics[width=0.5\textwidth]{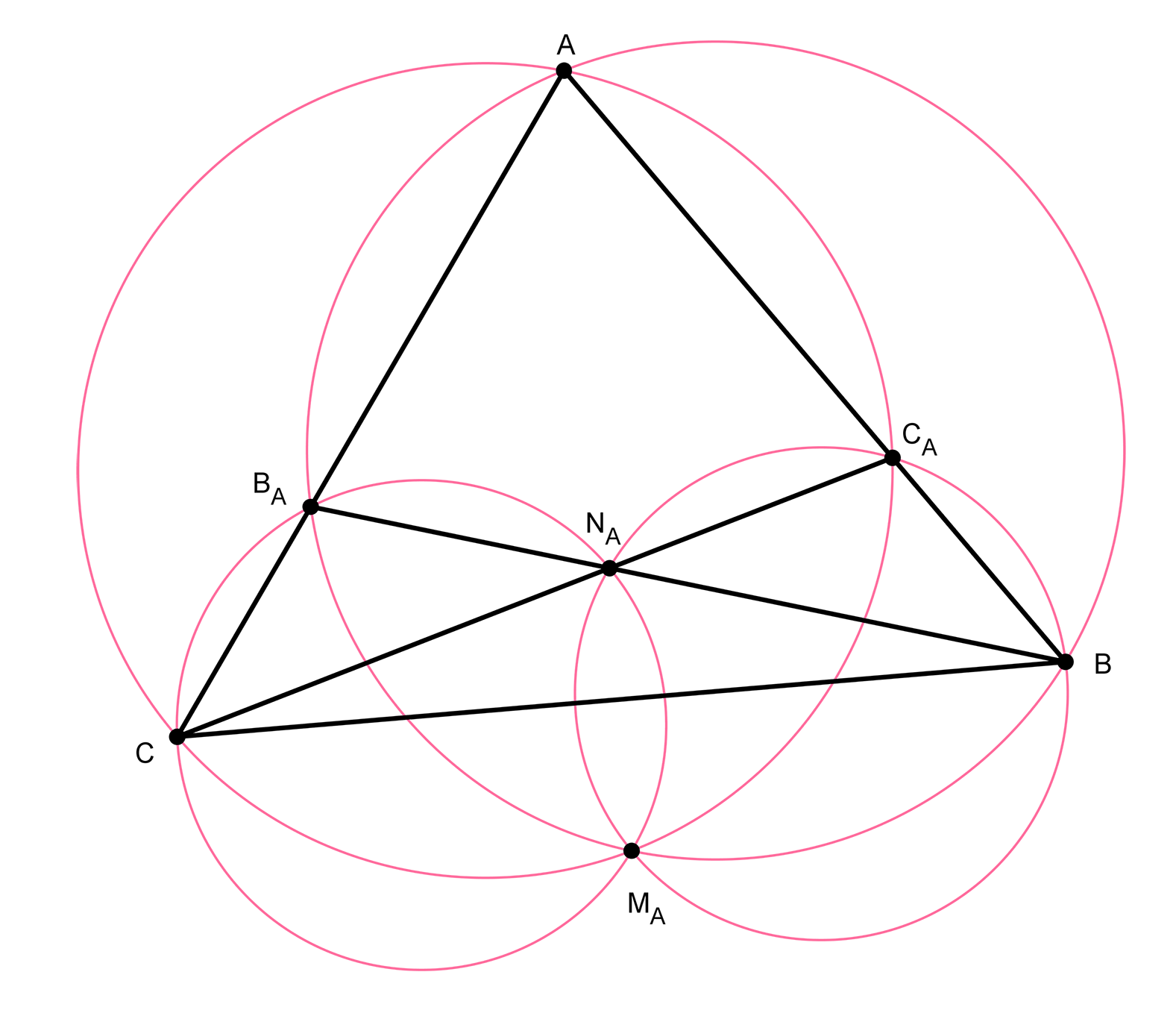}
			\newline
			Fig.3~Miquel-Steiner's point of triangle $\Delta{ABC}$ and cevians $BB_A, CC_A$
		\end{center}
		\begin{remark*}
			Cevians $BB_A, CC_A$ can be external. However we exclude situations $B_A=C$ or $C_A=B$ or $BB_A  \mid \mid CC_A$.
		\end{remark*}
		We can also reformulate previous lemma in the following way.
		\begin{lemma*}
			Miquel-Steiner's point $M_A$ belongs to the side $BC$ if and only if points $A, B_A, C_A, N_A$ are cyclic. 
		\end{lemma*}
		With this lemma we start research of Miquel-Steiner's point locus. 
	\section{Miquel-Steiner's point locus}
		Firstly we need to prove one lemma. 
		\begin{lemma}
			\label{lm:lm1}
			If cevians $BB_A, CC_A$ are parallel then circumcircles $\omega_{ACC_A}, \omega_{ABB_A}$ are tangent to each other (Fig.4).
		\end{lemma}
		\begin{center}
			\includegraphics[width=0.5\textwidth]{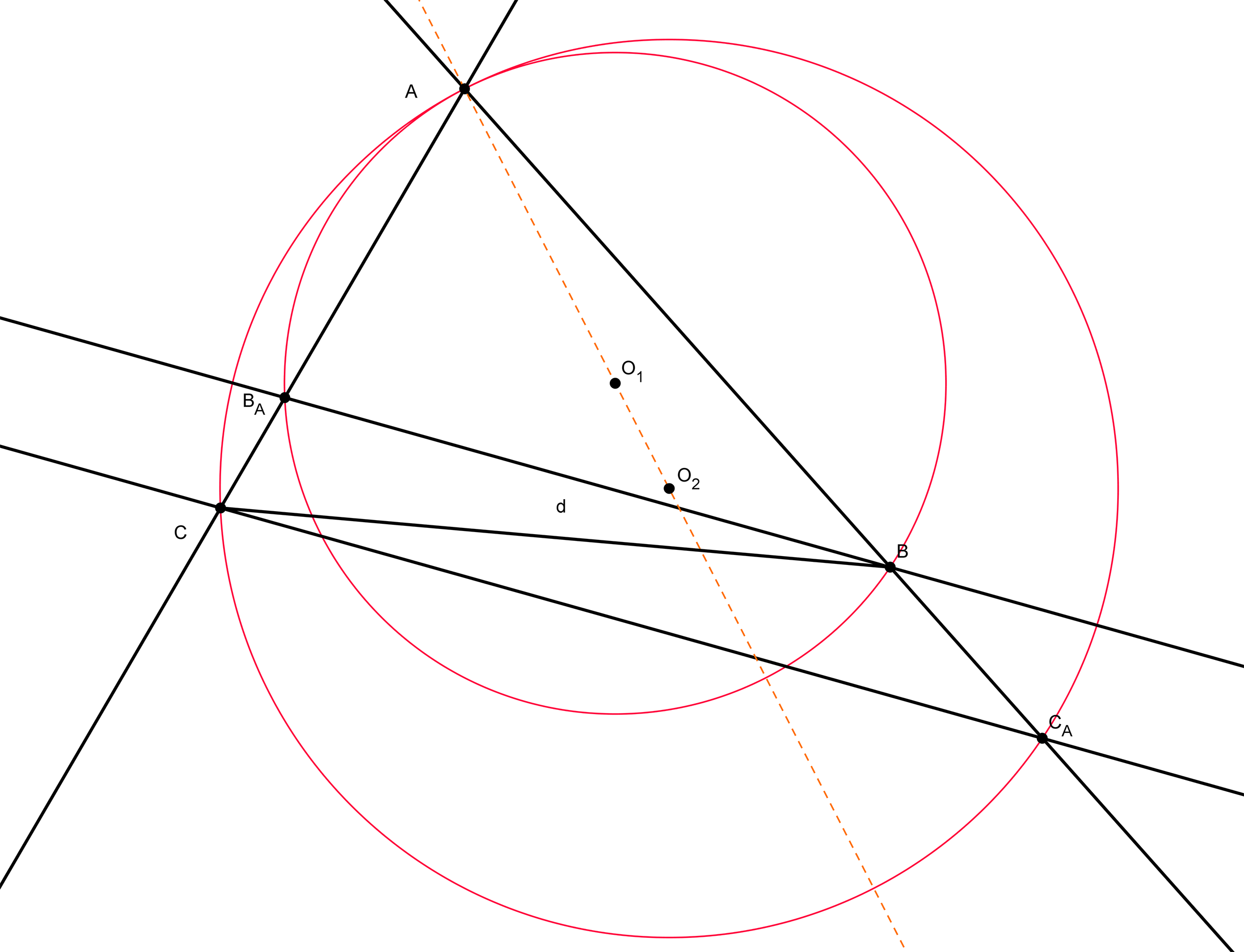}
			\newline
			Fig.4~Lemma 1
		\end{center}
		\begin{proof}
			Let us prove in a situation when $BB_A$ is internal. For other cases proof is analogical. 
			\newline
			$A, C, B_A$ are colinear, $A, B, C_A$ are colinear and $CC_A \mid \mid BB_A$. Therefore $\Delta{ABB_A}\sim\Delta{AC_AC}$. Let points $O_1, O_2$ be centres of $\omega_{ABB_A}, \omega_{ACC_A}$ respectively. Then $\angle{O_1AB_A}=\angle{O_2AC}$ (because of similarity). Thus $A, O_1, O_2$ are colinear. But $A$ is one of circumcircles intersection points. Thus $\omega_{ACC_A}, \omega_{ABB_A}$ are tangent to each other.
		\end{proof}
		\begin{statement}
			\label{st:st1}
			Let $\omega_{{ABC}}$ is circumcircle of $\Delta{ABC}$. 
			\newline
			1)~ Then for an arbitrary point $\forall{}M_A\not\in\omega_{{ABC}},M_A\not\in{AB},M_A\not\in{AC}:~\exists{!}$ pair $\left(BB_A, CC_A\right)$ such that $M_A$ is Miquel-Steiner's point for triangle ${\Delta{ABC}}$ with cevians $BB_A, CC_A$.
			\newline
			2)~$\forall$ pair $\left(BB_A, CC_A\right),B_A\neq{}C,C_A\neq{}B,BB_A  \nparallel CC_A:~\exists{!}M_A$ -- Miquel-Steiner's point.
		\end{statement}
		\begin{proof}
			1)~ $M_A$ is the candidate for Miquel-Steiner's point (Fig.3).
			Condition $M_A\not\in{AB},M_A\not\in{AC}$ makes possible to consider circumcicles $\omega_{ABM_A}, \omega_{ACM_A}$.
			Therefore $\exists{!}B_A=\omega_{ABM_A}\cap{AC}$ (intersection point different from $A$). Moreover $B_A\neq{C}$ because $M_A\not\in\omega_{ABC}$. Analogously $\exists{!}C_A=\omega_{ACM_A}\cap{AB}$ and $C_A\neq{B}$. Moreover $BB_A, CC_A$ cannot be parallel. If yes then by lemma \ref{lm:lm1} circumcircles $\omega_{ACC_A}, \omega_{ABB_A}$ are tangent to each other. It follows that $A=M_A$. But $M_A\not\in\omega_{ABC}$. Thus $BB_A \nparallel CC_A$. 
			\newline
			2)~ It is obvious. 
			\newline
			If $B_A\neq{}C$,$C_A\neq{}B$,$BB_A  \nparallel CC_A$ then points $A$,$C_A$,$N_A$,$B_A$,$B$,$C$ form complete quadrilateral. Therefore $M_A$ -- Miquel-Steiner's point exists and it is unique. 
			\qedhere
		\end{proof}
		\begin{corollary}
			\label{cor:cor1}
			There is one-to-one correspondence between admissible Miquel-Steiner's points and cevian pairs. 
		\end{corollary}
		Now we can see that locus in general case is entire plane without lines $AB, AC$ and circumcircle $\omega_{ABC}$. Therefore we need to add extra conditions.
		\begin{theorem}
			\label{th:th1}
			If cevians $BB_A, CC_A$ are internal then Miquel-Steiner's point locus is
			\begin{equation*}
				\Lambda_A=D_{ABC}\setminus \left(D_{A,B}\cap{}D_{A,C}\right),
			\end{equation*} 
			where $D_{ABC}$ is a disk bounded by circumcircle $\omega_{ABC}$, $D_{A,B}, D_{A,C}$ are disks bounded by $\omega_{A,B}, \omega_{A,C}$ respectively, where $\omega_{A,B}$ is tangent to $AB$ and  pass through $A$ and $C$, $\omega_{A,C}$ is tangent to $AC$ and  pass through $A$ and $B$ (Fig.5).
		\end{theorem}
		\begin{center}
			\includegraphics[width=0.5\textwidth]{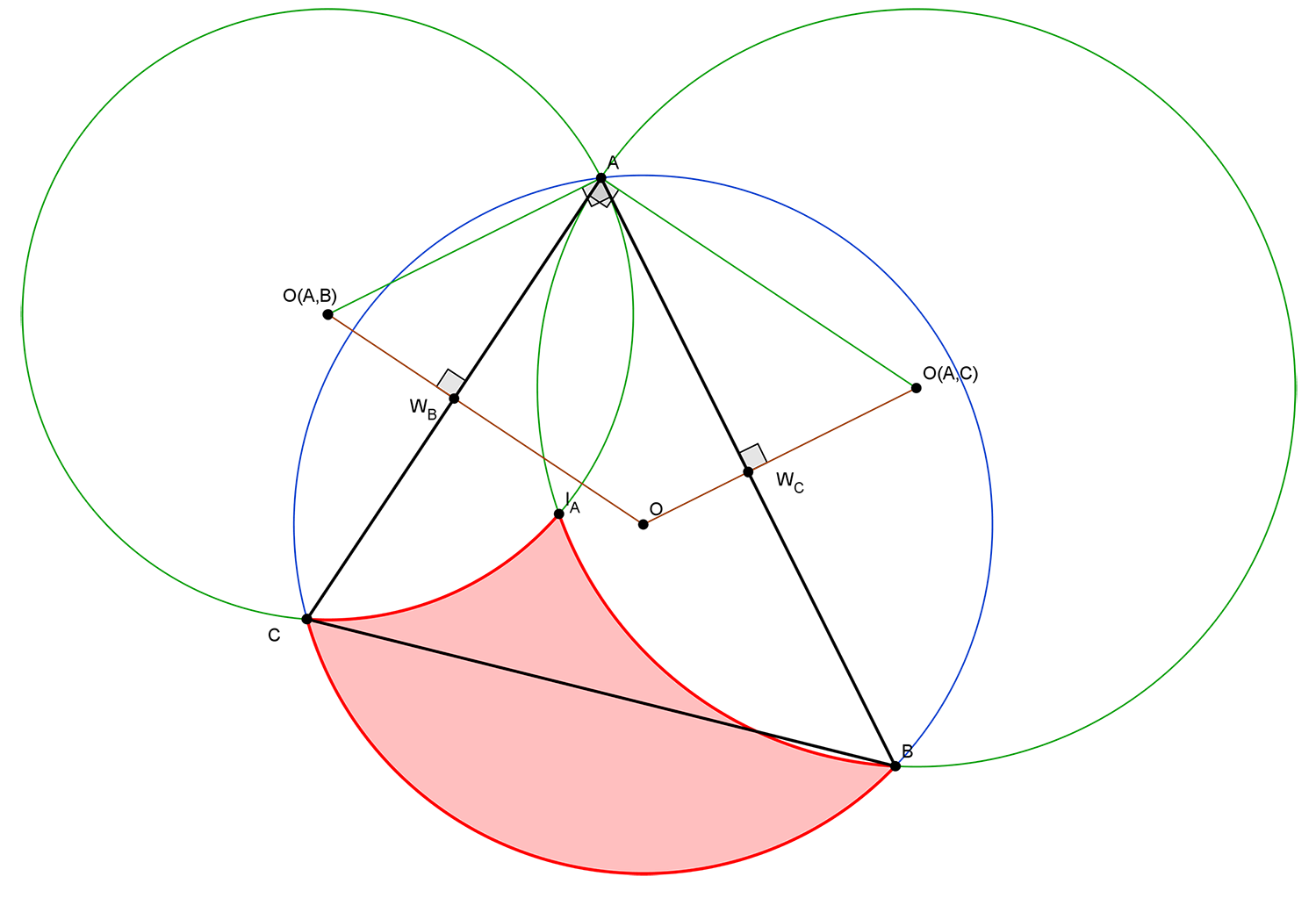}
			\newline
			Fig.5~Miquel-Steiner's point locus for internal cevians
		\end{center}
		\begin{proof}
			Let us prove in case when $\angle{A}$ is acute. For another case proof is analogical. Let us prove, that $CC_A$ internal only if $M_A\in{D_{ABC}\setminus{}D_{A,B}}$.
			We will do this step by step. 
			\newline
			Let semiplanes bounded by $AC$ will be $AC^{+}$ and $AC^{-}$ such that $B\in{AC^{+}}$. Analogously let semiplanes bounded by $AB$ will be $AB^{+}$ and $AB^{-}$ such that $C\in{AB^{+}}$
			\newline
			1)~Let us prove the following
			\begin{eqnarray*}
				& (a)~M_A\in{}D_{ABC}, M_A\in{AC^{+}} ~ \Rightarrow ~ \overrightarrow{BA} \upuparrows  \overrightarrow{BC_A} \\
				& (b)~M_A\not\in{}D_{ABC}, M_A\in{AC^{+}} ~ \Rightarrow ~ \overrightarrow{BA} \downarrow \uparrow  \overrightarrow{BC_A} \\
				& (c)~M_A\not\in{}D_{ABC}, M_A\in{AC^{-}} ~ \Rightarrow ~ \overrightarrow{BA} \upuparrows  \overrightarrow{BC_A} \\
				& (d)~M_A\in{}D_{ABC}, M_A\in{AC^{-}} ~ \Rightarrow ~ \overrightarrow{BA}  \downarrow \uparrow   \overrightarrow{BC_A}
			\end{eqnarray*}
			Let us prove $(a)$ proposition (Fig.6). Proof of $(b)$--$(d)$ will be analogical (considering instead of some angles their supplements). 
			\begin{center}
				\includegraphics[width=0.5\textwidth]{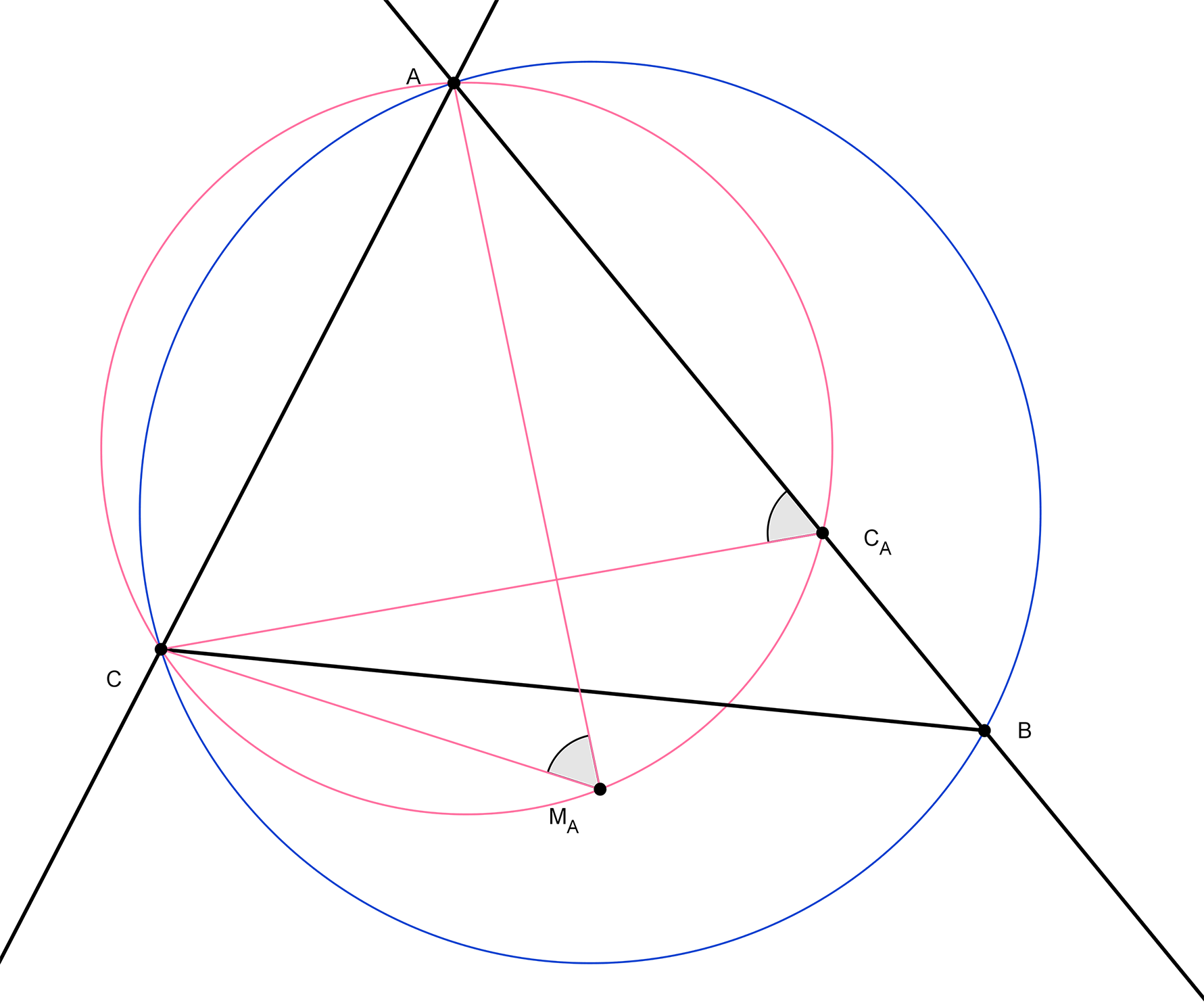}
				\newline
				Fig.6~$M_A\in{}D_{ABC}, M_A\in{AC^{+}}$
			\end{center}
			If $M_A\in{D_{ABC}, M_A\in{AC^{+}}}$ then $\angle{CM_AA}>\angle{CBA}$. If $C_A\in{AC^{-}}$ we will have needed result automatically. If not then $\angle{CM_AA}=\angle{CC_AA}$. Thus $\angle{CC_AA}>\angle{CBA}$. Finally 
			\begin{eqnarray*}
				&\angle{ACC_A}=180^\circ-\angle{BAC}-\angle{CC_AA}<180^\circ-\angle{BAC}-\angle{CBA}=\\
				&=\angle{ACB}~\Rightarrow ~ \overrightarrow{BA} \upuparrows  \overrightarrow{BC_A}
			\end{eqnarray*}
			2)~Let us prove the following
			\begin{eqnarray*}
				& (a)~M_A\in{}D_{A,B}, M_A\in{AC^{+}} ~ \Rightarrow ~ \overrightarrow{AB} \downarrow \uparrow  \overrightarrow{AC_A} \\
				& (b)~M_A\not\in{}D_{A,B}, M_A\in{AC^{+}} ~ \Rightarrow ~ \overrightarrow{AB} \upuparrows  \overrightarrow{AC_A} \\
				& (c)~M_A\not\in{}D_{A,B}, M_A\in{AC^{-}} ~ \Rightarrow ~ \overrightarrow{AB} \downarrow \uparrow  \overrightarrow{AC_A} \\
				& (d)~M_A\in{}D_{A,B}, M_A\in{AC^{-}} ~ \Rightarrow ~ \overrightarrow{AB}  \upuparrows   \overrightarrow{AC_A}
			\end{eqnarray*}
			Let us prove $(a)$ proposition (Fig.7). Proof of $(b)$--$(d)$ will be analogical (considering instead of some angles their supplements). 
			\begin{center}
				\includegraphics[width=0.5\textwidth]{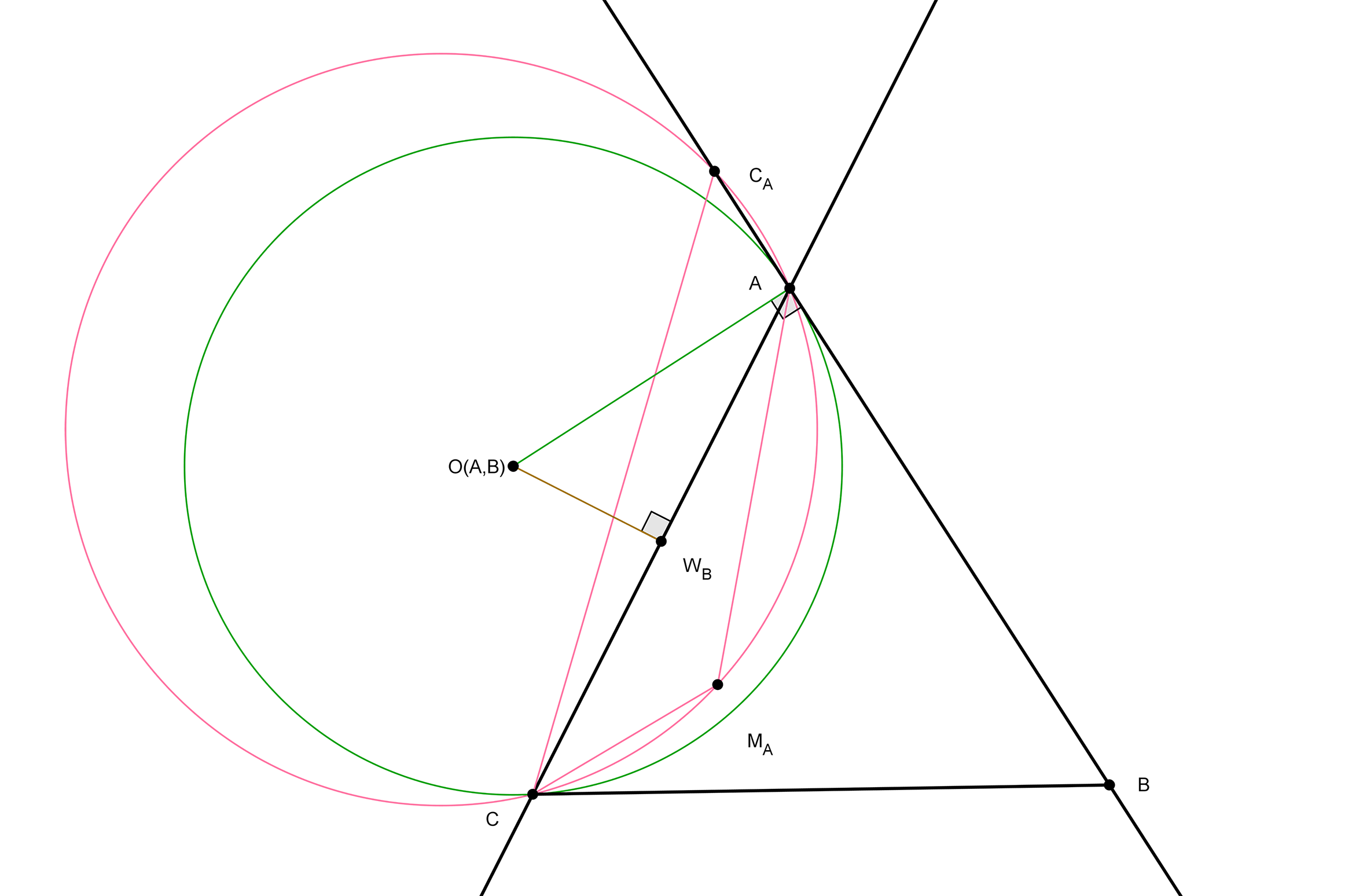}
				\newline
				Fig.7~$M_A\in{}D_{A,B}, M_A\in{AC^{+}}$
			\end{center}
			If $M_A\in{}D_{A,B}, M_A\in{AC^{+}}$ then $\angle{CM_AA}>180^{\circ}-\angle{BAC}$. Assuming that $C_A \in {AC^{+}}$ we will have $\angle{CC_AA}=\angle{CM_AA}$ and $\angle{ACC_A}=180^{\circ}-\angle{BAC}-\angle{CC_AA}$. It follows that $\angle{ACC_A}<0$. Contradiction. Therefore $C_A \in {AC^{-}}$ and we have proved needed. 
			\newline
			Therefore $CC_A$ internal only if $M_A\in{D_{ABC}\setminus{}D_{A,B}}$. Analogously $BB_A$ internal only if $M_A\in{D_{ABC}\setminus{}D_{A,C}}$.
			\newline
			Thus $BB_A, CC_A$ are internal only if $M_A\in{\Lambda_A}$. With corollary \ref{cor:cor1} we have converse proposition. Therefore $\Lambda_A$ is locus. 
			\qedhere
		\end{proof}
		\begin{remark}
			\label{rm:rm1}
			Circles $\omega_{ABC}, \omega_{A,B}, \omega_{A,C}$ are actually limiting circles
				\begin{eqnarray*}
					& {\omega_{AC_AC}}\to\omega_{ABC},~where~{C_A \to B, C_A\in{AB}} \\
					& {\omega_{AC_AC}}\to\omega_{A,B},~where~{C_A \to A, C_A\in{AB}} \\
					& {\omega_{ABB_A}}\to\omega_{ABC},~where~{B_A \to C, B_A\in{AC}} \\
					& {\omega_{ABB_A}}\to\omega_{A,C},~where~{B_A \to A, B_A\in{AC}}
				\end{eqnarray*}
			Therefore it is natural to consider them.
		\end{remark}
		\begin{definition}
			\label{def:def1}
			Let circles $\omega_{A,B}, \omega_{A,C}$ be called Miquel-Steiner's auxiliary circles of $\Delta{ABC}$ from $A$. Let centres $O_{A,B}, O_{A,C}$ of circles $\omega_{A,B}, \omega_{A,C}$ be called auxiliary Miquel-Steiner's centres of $\Delta{ABC}$ from $A$. Let $I_A=\omega_{A,B}\cap\omega_{A,C}\neq{A}$ be called main Miquel-Steiner's centre of $\Delta{ABC}$ from $A$. Let line $AI_A$ be called Miquel-Steiner's axis of $\Delta{ABC}$ from $A$ (Fig.8)
		\end{definition}
		\begin{center}
			\includegraphics[width=0.5\textwidth]{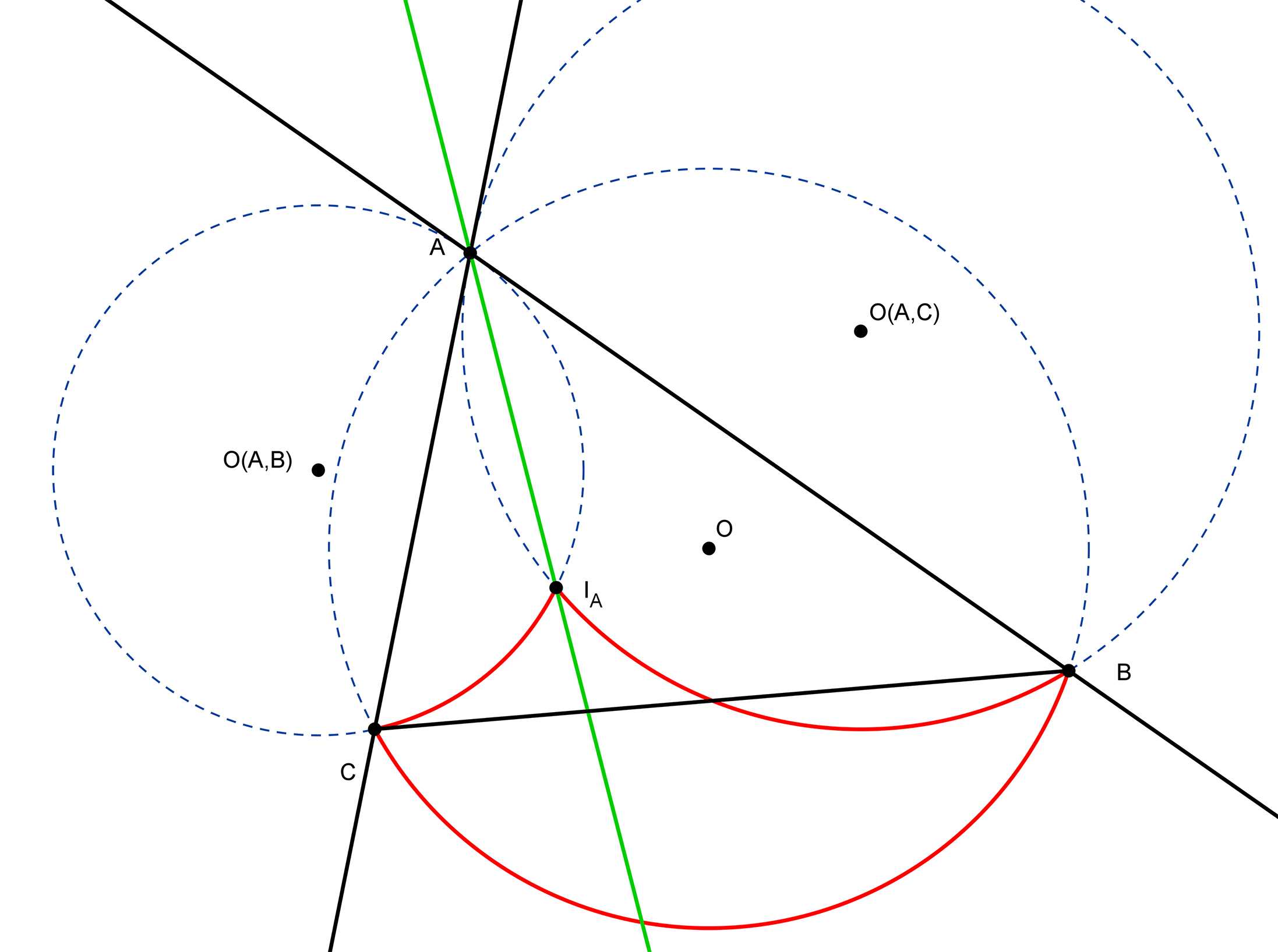}
			\newline
			Fig.8~Miquel-Steiner's main and auxiliary centres and Miquel-Steiner's axis
		\end{center}
	\section{Miquel-Steiner's point locus and the Brocard's circle}
		It is important to notice, that the following results (theorems \ref{th:th2}--\ref{th:th4} and lemma \ref{lm:lm2}) were obtained and proved in 2012 as a part of Miquel-Steiner's point locus study (the All-Ukrainian Research Paper Defense Competition Among Student-members of Minor Academy of Sciences of Ukraine). However it was found that these results were not original. They had been obtained earlier by french mathematician Henri Brocard. 
		\newline
		Miquel-Steiner's auxiliary circles are also known as conjugate circles. ~[1, p.190] 
		\begin{theorem}[Intermediate results of Brocard contributions*]
			\label{th:th2}
			Let point $O$ be circumcircle centre of $\Delta{ABC}$. Then
			\newline
			1)~Centres $O_{A,B}, O_{A,C}, I_A, O$ are cyclic
			\newline
			2)~$O_{A,B}O_{A,C} \mid \mid OI_A \perp AI_A$ (Fig.9) 
		\end{theorem}
		\begin{center}
			\includegraphics[width=0.5\textwidth]{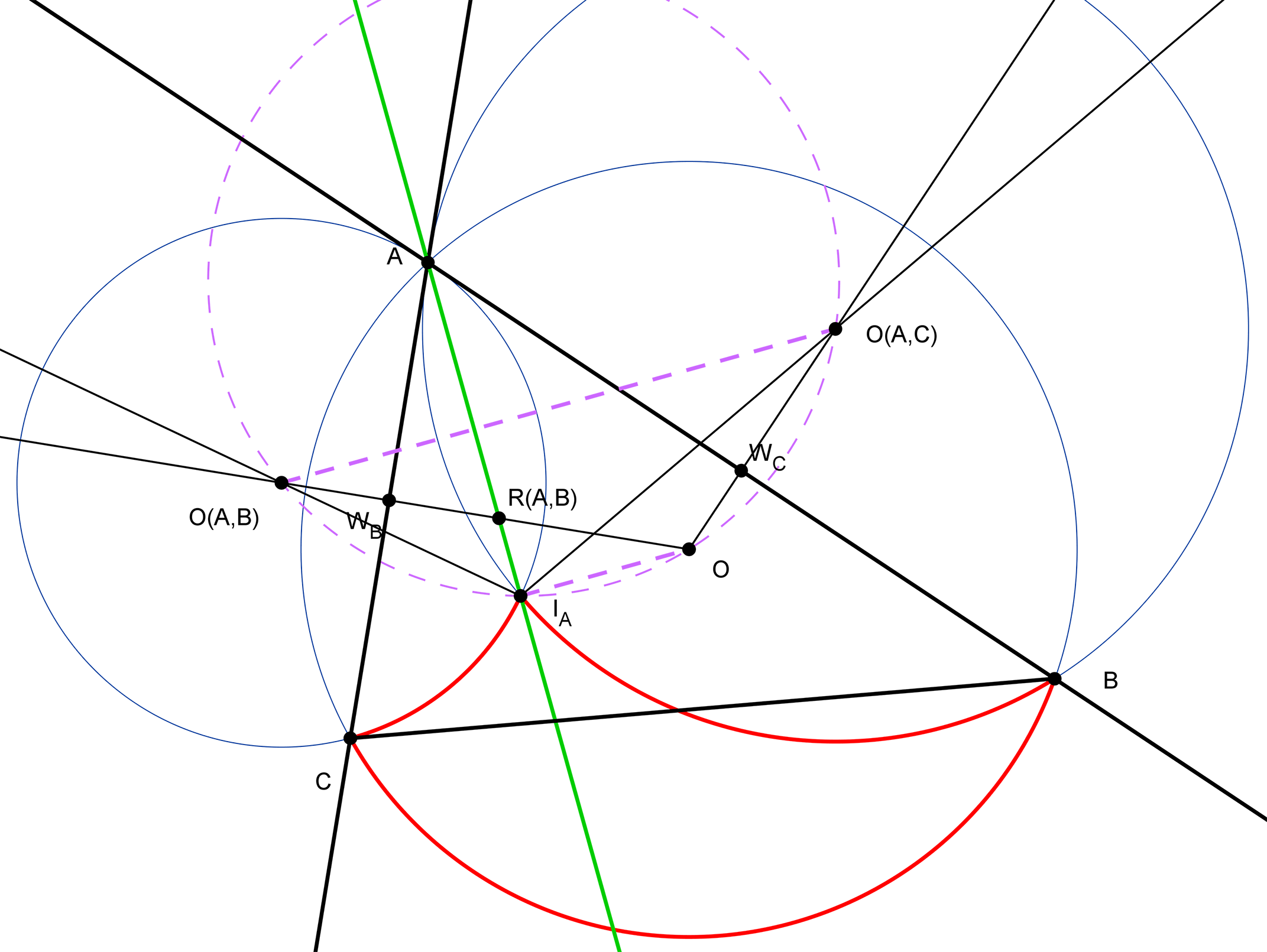}
			\newline
			Fig.9~Theorem 2
		\end{center}
		Following lemma is limiting case of the Miquel's theorem ~[2, p.61]. 
		\begin{lemma}[Brocard's points*]
			\label{lm:lm2}
			Auxiliary Miquel-Steiner's circles $\omega_{A,C}$, $\omega_{B,A}$, $\omega_{C,B}$ have similar point $G$ (Fig.10)
		\end{lemma}
		\begin{center}
			\includegraphics[width=0.5\textwidth]{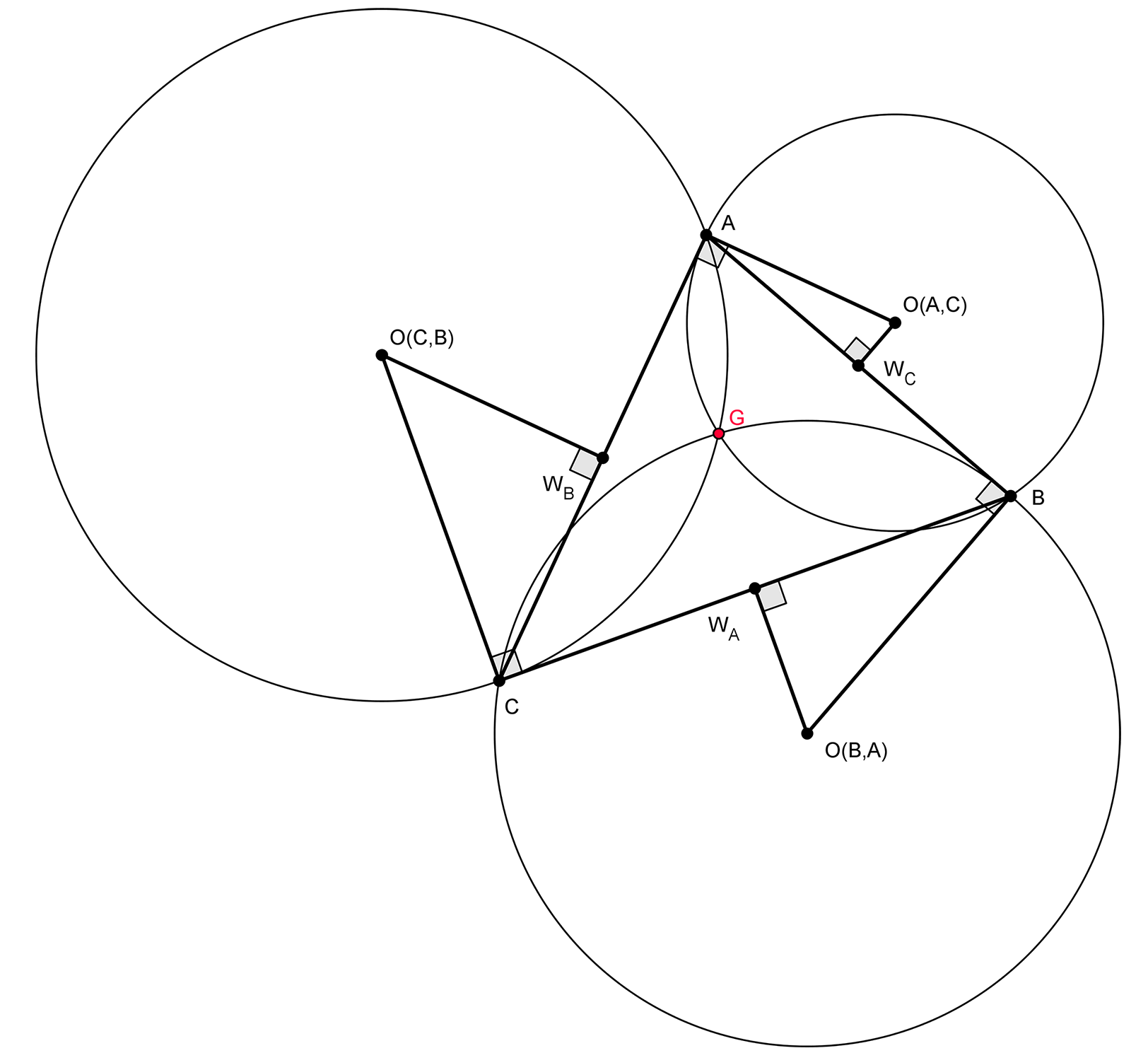}
			\newline
			Fig.10~Lemma 2
		\end{center}
		\begin{remark}
			\label{rm:rm2}
			Triangle $\Delta{ABC}$ has two points from lemma \ref{lm:lm2}. They are also known as Brocard's points ~[3, pp.98-124].  
		\end{remark}
		Now we can consider theorem. 
		\begin{theorem}[Brocard's circle and Brocard's second triangle*]
			\label{th:th3}
			$I_A, I_B, I_C$--three Miquel-Steiner's main centres of $\Delta{ABC}$ are cyclic with $O$. (Fig.11)
		\end{theorem}
		\begin{center}
			\includegraphics[width=0.5\textwidth]{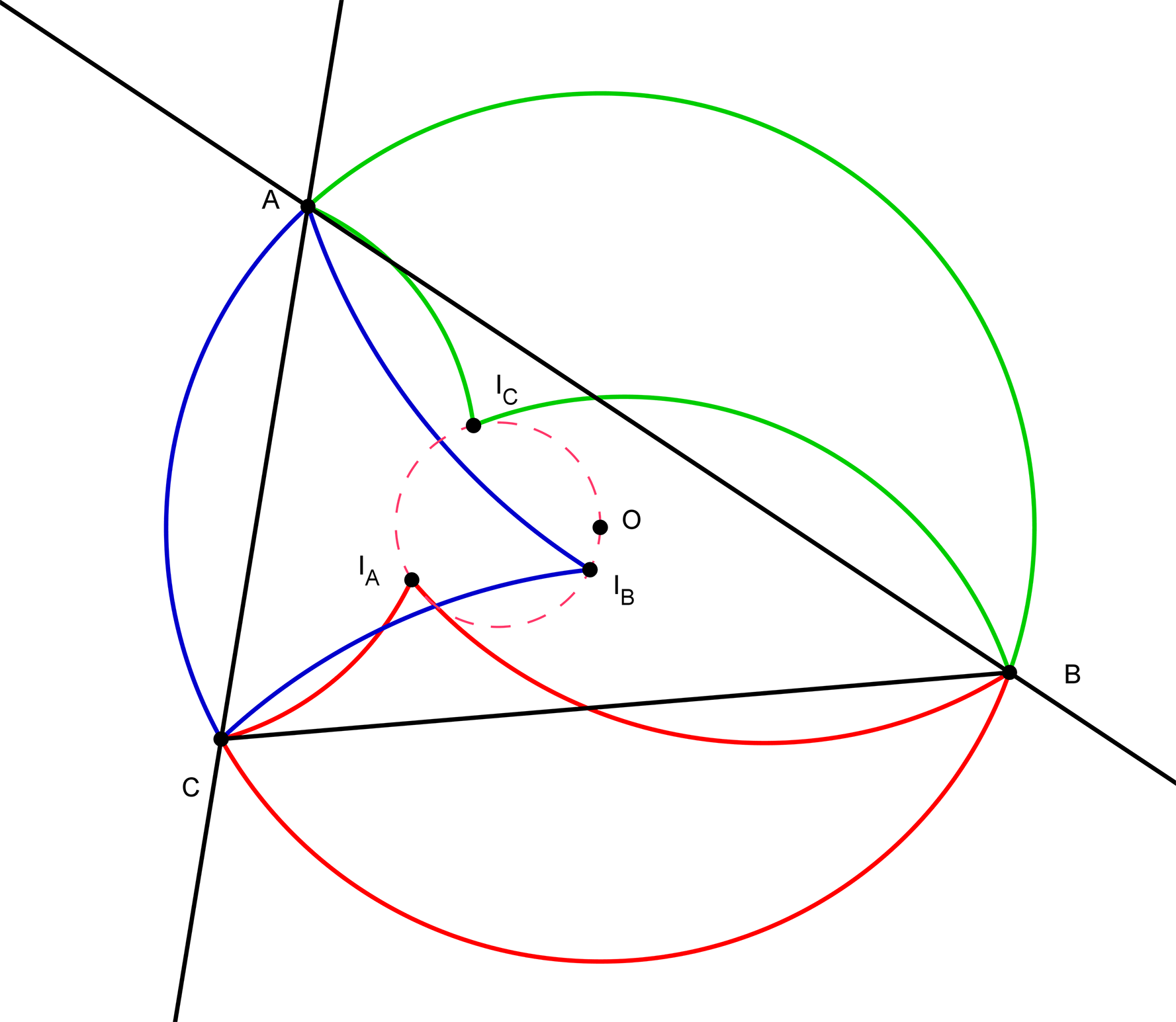}
			\newline
			Fig.11~Theorem 3
		\end{center}
		\begin{remark}
			\label{rm:rm3}
			Points $I_A, I_B, I_C$ from theorem \ref{th:th3} are also known as vertices of the Brocard's second triangle ~[3, pp.110-118]. Circle $\omega_{I_AI_BI_CO}$ from theorem \ref{th:th3} is also known as the Brocard's circle ~[3, pp.106-110].
		\end{remark}
		Miquel-Steiner's axises are not simple lines. 
		\begin{theorem}[Intermediate results of Brocard contributions*]
			\label{th:th4}
			Miquel-Steiner's axises are symmedians. (Fig.12)
		\end{theorem}
		\begin{center}
			\includegraphics[width=0.5\textwidth]{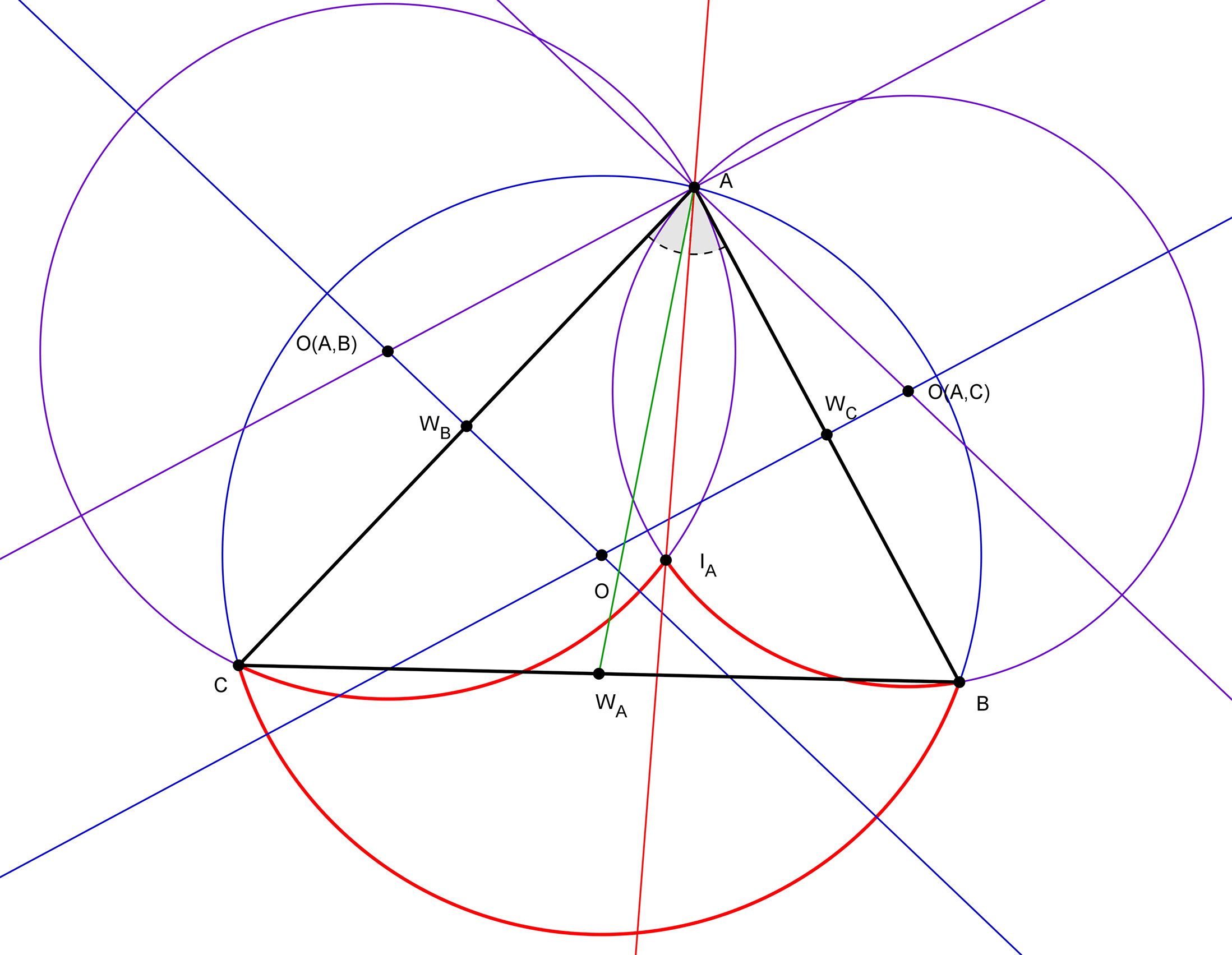}
			\newline
			Fig.12~Theorem 4
		\end{center}
	\section{Some specific locuses}
		Here we consider cases when cevians has some specific properties.
		\begin{statement}
			\label{st:st2}
			Let $M$ is intersection point of tangents to $\omega_{ABC}$ in points $B,C$. Let $\omega_{tan}$ is circle of $B,C$ with centre $M$. 
			\newline
			Then cevians $BB_A, CC_A$ are perpendicular if and only if Miquel-Steiner's point $M_A$ belongs to $\omega_{tan}$. (Fig.13) 
		\end{statement}
		\begin{center}
			\includegraphics[width=0.5\textwidth]{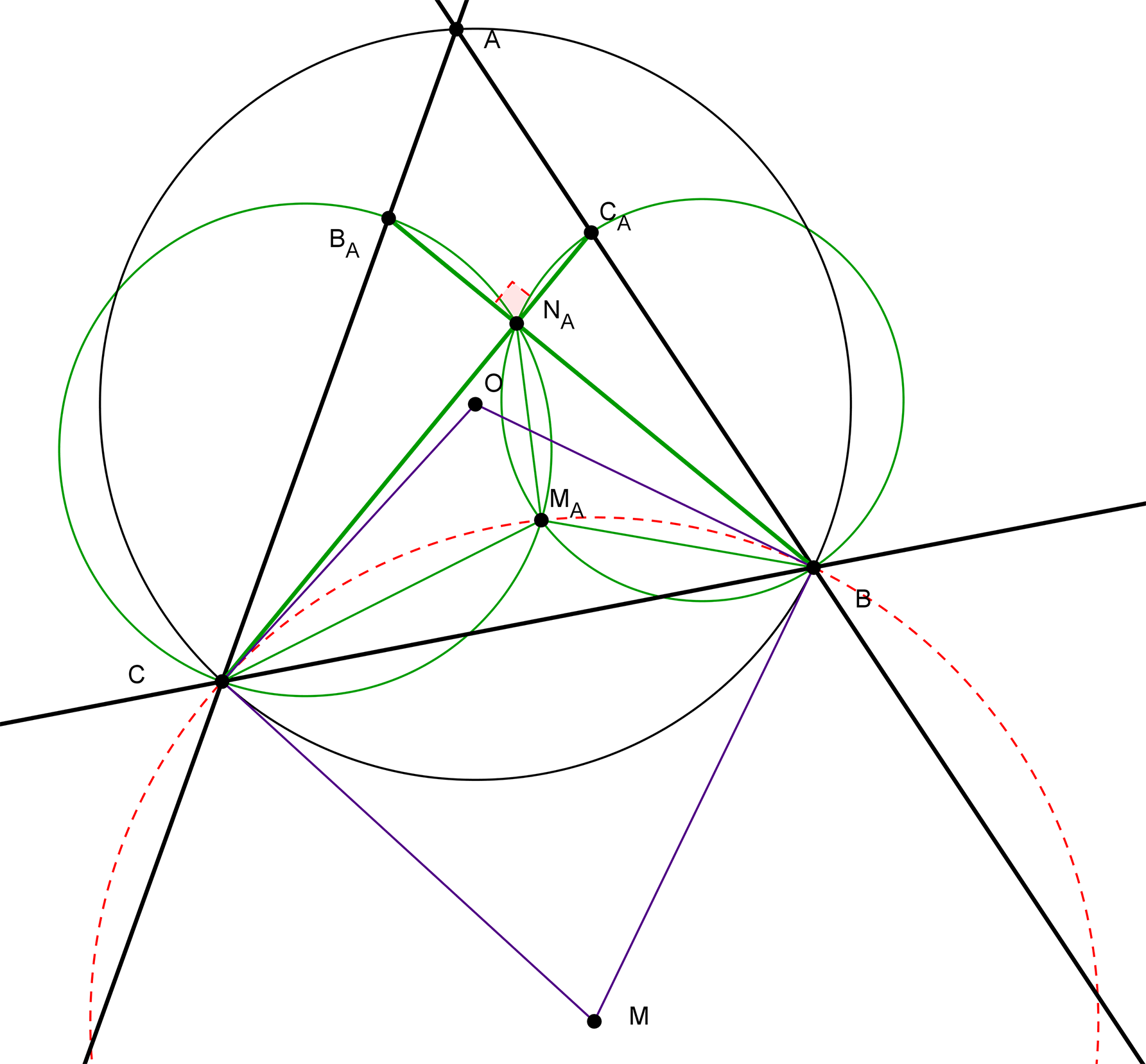}
			\newline
			Fig.13~Statement 2
		\end{center}
		\begin{proof}
			Let us proof in case when $\angle{A}$ is acute. For another case proof is analogical. Finally, let us proof when cevians are internal. For external or mixed cevians proof is analogical.
			\newline
			Obvious that 
			\begin{equation*}
				\angle{CMB}=180^{\circ}-\angle{BOC}=180^{\circ}-2\angle{BAC} 
			\end{equation*}
			With this we have 
			\begin{equation}
				\label{eq:eq11}
				M_A\in{\omega_{tan}}~\Leftrightarrow~\angle{BM_AC}=180^{\circ}-\frac{1}{2}\angle{CMB}=90^{\circ}+\angle{BAC} 
			\end{equation}
			Also because we have
			\begin{eqnarray*}
				&\angle{BM_AC}=360^{\circ}-\angle{CM_AN_A}-\angle{N_AM_AB}=\\
				&=\angle{N_AB_AC}+\angle{BC_AN_A}=\angle{B_AN_AC_A}+\angle{BAC}
			\end{eqnarray*}
			It follows that
			\begin{equation}
				\label{eq:eq12}
				\angle{B_AN_AC_A}=90^{\circ}~\Leftrightarrow~\angle{BM_AC}=90^{\circ}+\angle{BAC} 
			\end{equation}
			With (\ref{eq:eq11})--(\ref{eq:eq12}) we have needed.
			\qedhere
		\end{proof}
		\begin{statement}
			\label{st:st3}
			$B_AC_A \mid\mid BC$ if and only if Miquel-Steiner's point $M_A$ belongs to symmedian from $A$ (excluding $A$). (Fig.14)
		\end{statement}
		\begin{center}
			\includegraphics[width=0.5\textwidth]{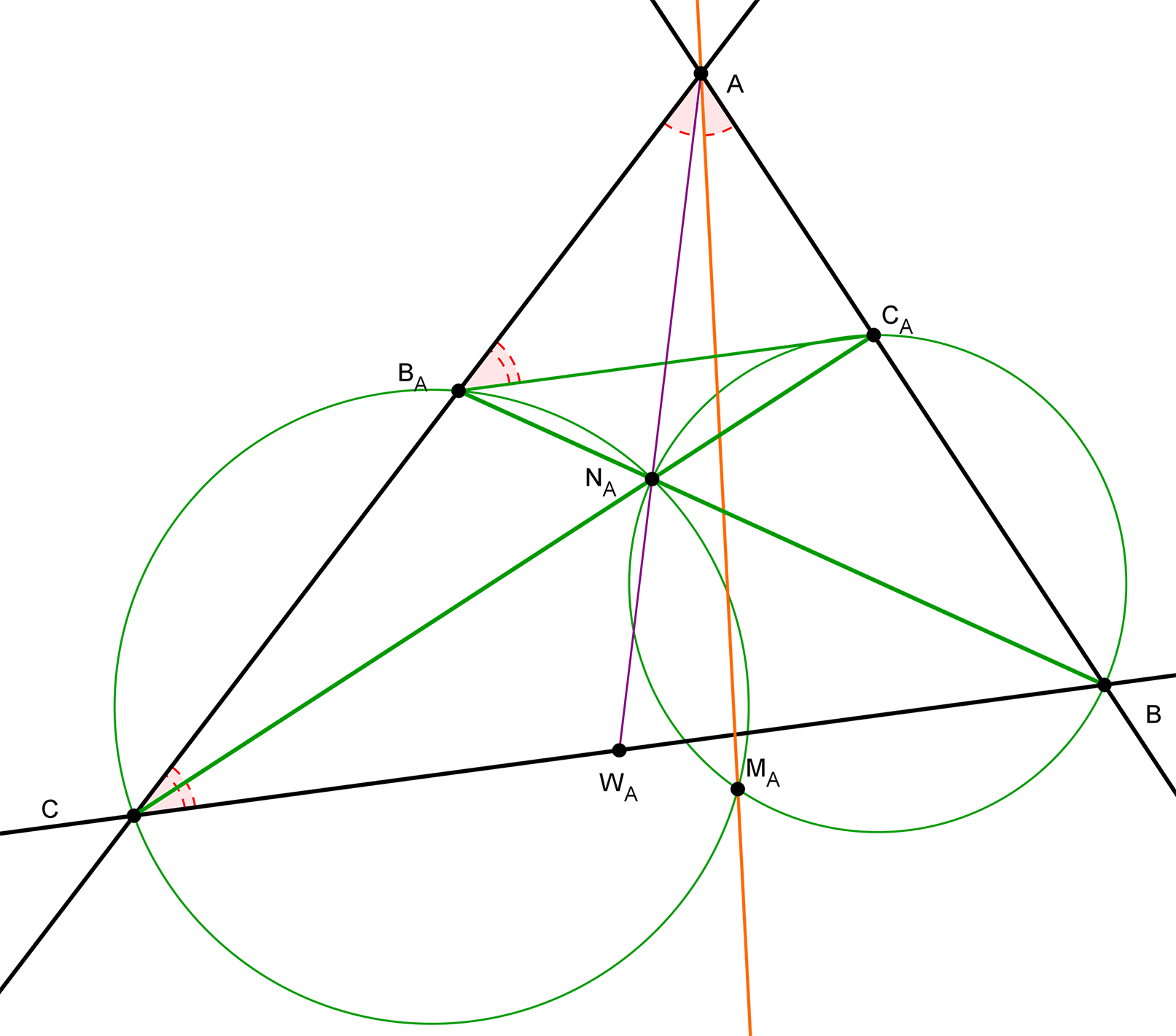}
			\newline
			Fig.14~Statement 3
		\end{center}
		\begin{proof}
			Let us proof in case when $\angle{A}$ is acute. For another case proof is analogical. Proof has four steps.
			\newline
			1)~Let us proof that if $B_AC_A \mid\mid BC$ and $B'_AC'_A \mid\mid BC$ then Miquel-Steiner's points $M_A, M'_A$ are colinear with $A$. Let this line be $l_{par}$ (Fig.15)
			\begin{center}
				\includegraphics[width=0.5\textwidth]{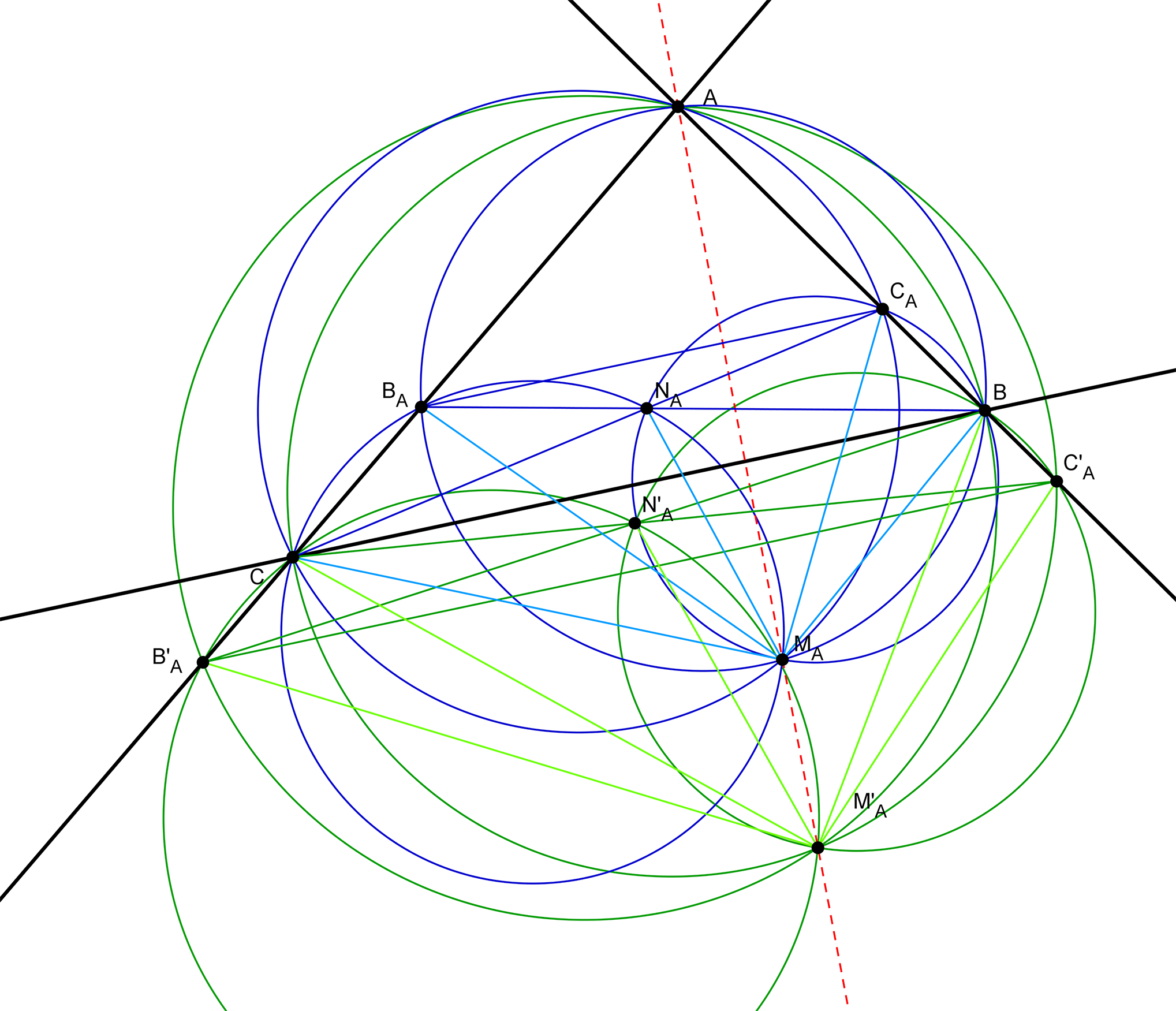}
				\newline
				Fig.15~Statement 3, proof, step 1
			\end{center}
			Let us notice following angles equalities
			\begin{eqnarray*}
				& \angle{CM_AC_A}=180^{\circ}-\angle{BAC}=\angle{CM'_AC'_A} \\
				& \angle{M_AB_AC}=\angle{M_AN_AC}=\angle{M_ABA} \\
				& \angle{B_ACM_A}=\angle{BN_AM_A}=\angle{BC_AM_A} \\
				& \angle{CB'_AM'_A}=\angle{C'_AN'_AM'_A}=\angle{C'_ABM'_A} \\
				& \angle{M'_ACB'_A}=\angle{M'_AN'_AB'_A}=\angle{M'_AC'_AB}
			\end{eqnarray*}
			It follows that
			\begin{equation*}
				\Delta{CB_AM_A}\sim\Delta{C_ABM_A},~\Delta{CB'_AM'_A}\sim\Delta{C'_ABM'_A}
			\end{equation*}
			Thus 
			\begin{eqnarray*}
				& \frac{\left|CM_A\right|}{\left|C_AM_A\right|}=\frac{\left|CB_A\right|}{\left|C_AB\right|}=\frac{\left|AC\right|}{\left|AB\right|} \\
				& \frac{\left|CM'_A\right|}{\left|C'_AM'_A\right|}=\frac{\left|CB'_A\right|}{\left|BC'_A\right|}=\frac{\left|AC\right|}{\left|AB\right|}
			\end{eqnarray*}
			It follows that
			\begin{equation*}
				\Delta{CM_AC_A}\sim\Delta{CM'_AC'_A}
			\end{equation*}
			Thus $\angle{CC'_AM'_A}=\angle{CC_AM_A}$. It follows that
			\begin{equation*}
				\angle{CAM'_A}=\angle{CC'_AM'_A}=\angle{CC_AM_A}=\angle{CAM_A}
			\end{equation*}
			Therefore $M_A, M'_A, A$ are colinear.
			\newline
			2)~Let $M$ is intersection point of tangents to $\omega_{ABC}$ in points $B,C$ (as in statement \ref{st:st2}). Let $M$ is considered as Miquel-Steiner's point with respective cevians $BD, CE$. Then $D, E, M$ are colinear and $ED \mid\mid BC$. (Fig.16)
			\begin{center}
				\includegraphics[width=0.5\textwidth]{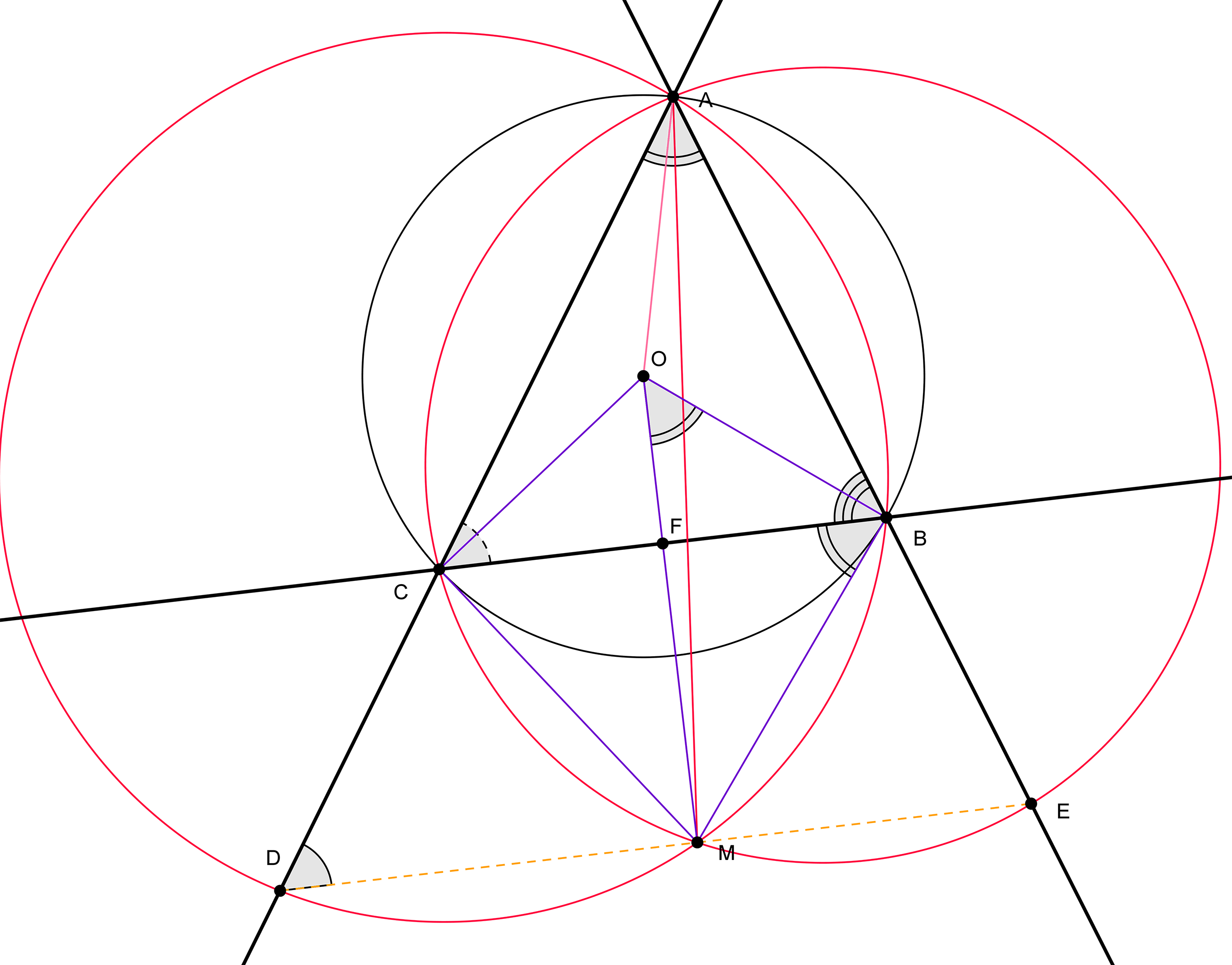}
				\newline
				Fig.16~Statement 3, proof, step 2
			\end{center}
			It is obvious that $\Delta{OCM}=\Delta{OBM}$. It follows that $BOCM$ is deltoid and if $F=BC\cap{OM}$ then $\angle{OFB}=90^{\circ}$.
			It follows that
			\begin{equation*}
				\angle{CBM}=90^{\circ}-\angle{BMF}=\angle{FOB}=\frac{1}{2}\angle{COB}=\angle{CAB}
			\end{equation*}
			But 
			\begin{equation*}
				\angle{ADM}=180^{\circ}-\angle{ABM}=180^{\circ}-\angle{ABC}-\angle{CBM}
			\end{equation*}
			Finally
			\begin{equation*}
				\angle{ADM}=\angle{ACB}
			\end{equation*}
			It follows that $MD\mid\mid{BC}$. Analogously $ME\mid\mid BC$. Therefore $M, E, D$ are colinear and $ED\mid\mid BC$.
			\newline
			3)~With step 2 we have that if $M$ is Miquel-Steiner's point then cevian bases $D, E$ such that $ED \mid\mid BC$. With step 1 it means that $M\in{l_{par}}$. But $A\in{l_{par}}$. Therefore $l_{par}=AM$. 
			\newline
			From the other hand with one narrowly known symmedian property [1, p.150] $AM$ is symmedian from $A$ ($MB, MC$ are also called external symmedians).
			\newline 
			Therefore we have proved that if $B_AC_A \mid\mid BC$ then Miquel-Steiner's point belongs to $l_{par}=AM$ which is symmedian from $A$. 
			\newline
			4)~Let us prove converse proposition. Or if Miquel-Steiner's point belongs to symmedian then $B_AC_A \mid\mid BC$. (Fig.17)
			\begin{center}
				\includegraphics[width=0.5\textwidth]{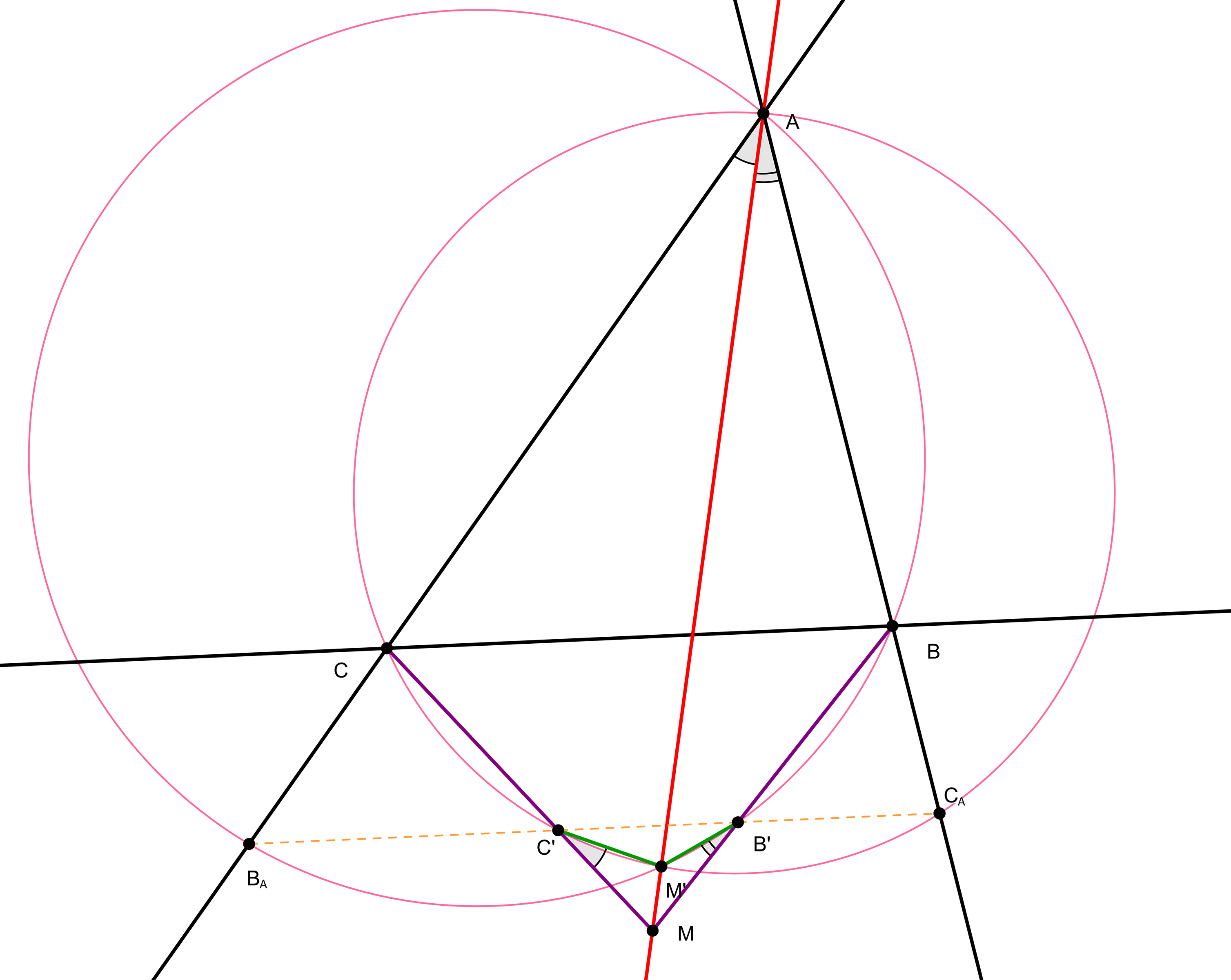}
				\newline
				Fig.17~Statement 3, proof, step 4
			\end{center}
			Let $M'$ -- another Miquel-Steiner's point on symmedian $AM$. Let $B'=\omega_{ABM'} \cap MB$, $C'=\omega_{ACM'} \cap MC$. Therefore $A,C,C',M'$ are cyclic. It follows that $\angle{M'C'M}=\angle{MAC}$. Analogously $\angle{MB'M'}=\angle{BAM}$. Therefore we have
			\begin{equation*}
				\Delta{M'C'M}\sim\Delta{CAM},~\Delta{M'B'M}\sim\Delta{BAM}
			\end{equation*}
			With this similarity we have
			\begin{equation*}
				\frac{\left|MC'\right|}{\left|AM\right|}=\frac{\left|MM'\right|}{\left|MC\right|},~\frac{\left|MB'\right|}{\left|AM\right|}=\frac{\left|MM'\right|}{\left|MB\right|}
			\end{equation*}
			But $\left|MB\right|=\left|MC\right|$. Therefore $\left|MB'\right|=\left|MC'\right|$. It follows that $ B'C'\mid\mid BC$. From step 2 we have $\angle{MBC}=\angle{BAC}$. Therefore
			\begin{equation*}
				\angle{BB'C'}= 180^{\circ}-\angle{MBC}=180^{\circ}-\angle{BAC}
			\end{equation*}
			From the other hand if $B_A=\omega_{ABM'}\cap{AC}$ then $\angle{BB'B_A}=180^{\circ}-\angle{BAC}$. Thus $\angle{BB'B_A}=\angle{BB'C'}$ and $B_A,B',C'$ are colinear. Analogously $C_A,B',C'$.
			\newline
			Therefore $B_AC_A\mid\mid BC$. 
			\qedhere
		\end{proof}
		\begin{remark}
			\label{rm:rm4}
			With one trapezoid property ~[4, p.23-35] (also known as Steiner's line theorem) we can reformulate statement \ref{st:st3}. \newline
			Miquel-Steiner's point $M_A\neq{A}$ belongs to symmedian from $A$ if and only if cevians intersection point $N_A\neq{A}$ belongs to median from $A$. (Fig.18)
		\end{remark}
		\begin{center}
			\includegraphics[width=0.5\textwidth]{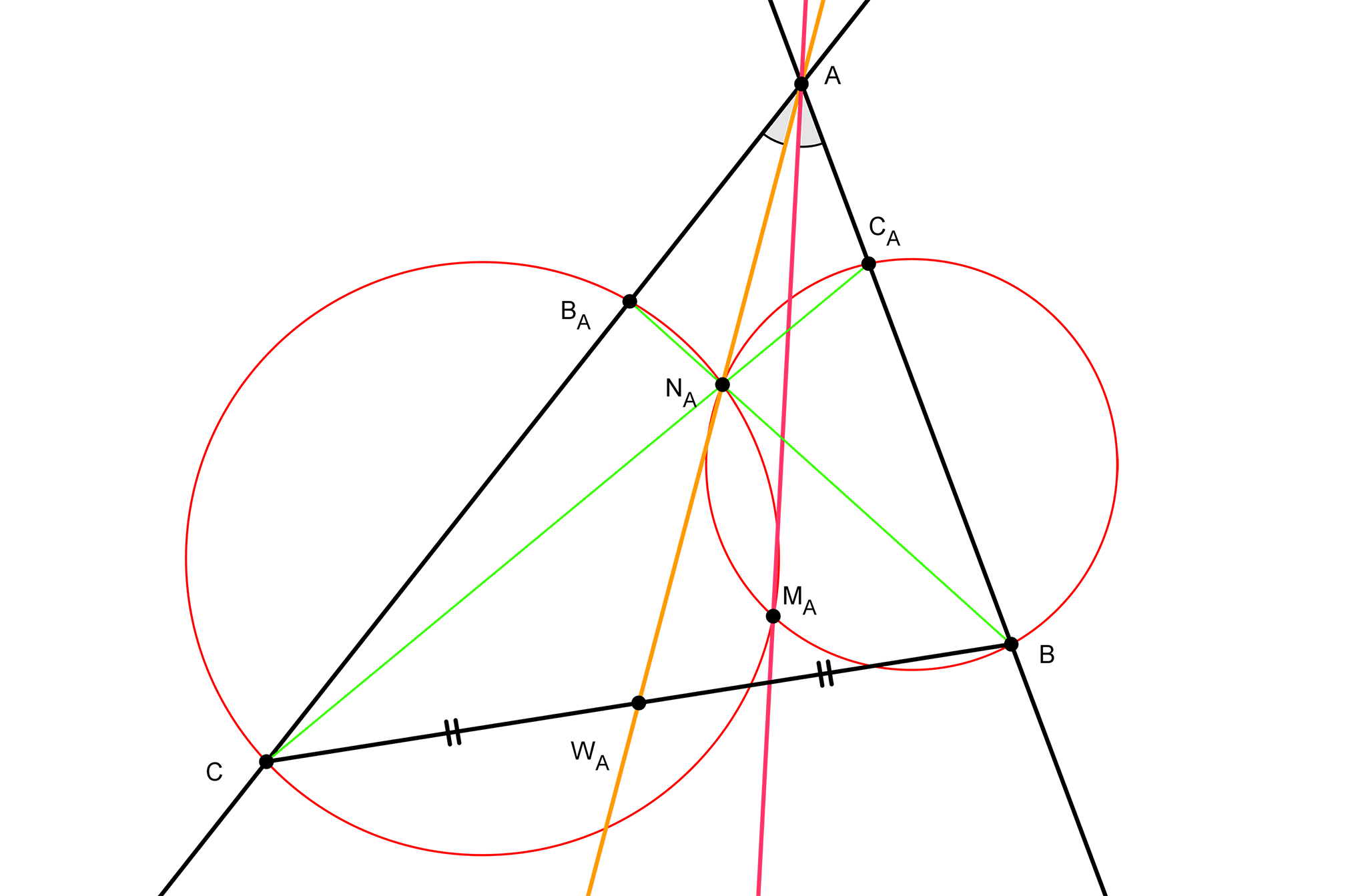}
			\newline
			Fig.18~Miquel's Steiner point belongs to symmedian
		\end{center}
		Median is isogonal to symmedian. It is logical to generalize previous statement in a case, when the Miquel-Steiner's point belongs to an arbitrary line from $A$. 
		\begin{statement}
			\label{st:st4}
			If Miquel-Steiner's point $M_A\neq{A}$ belongs to fixed line $l$ and $A\in{l}$ then respective cevians intersection point $N_A$ belongs to line which is parallel to isogonal of $l$. (Fig.19)
		\end{statement}
		\begin{center}
			\includegraphics[width=0.5\textwidth]{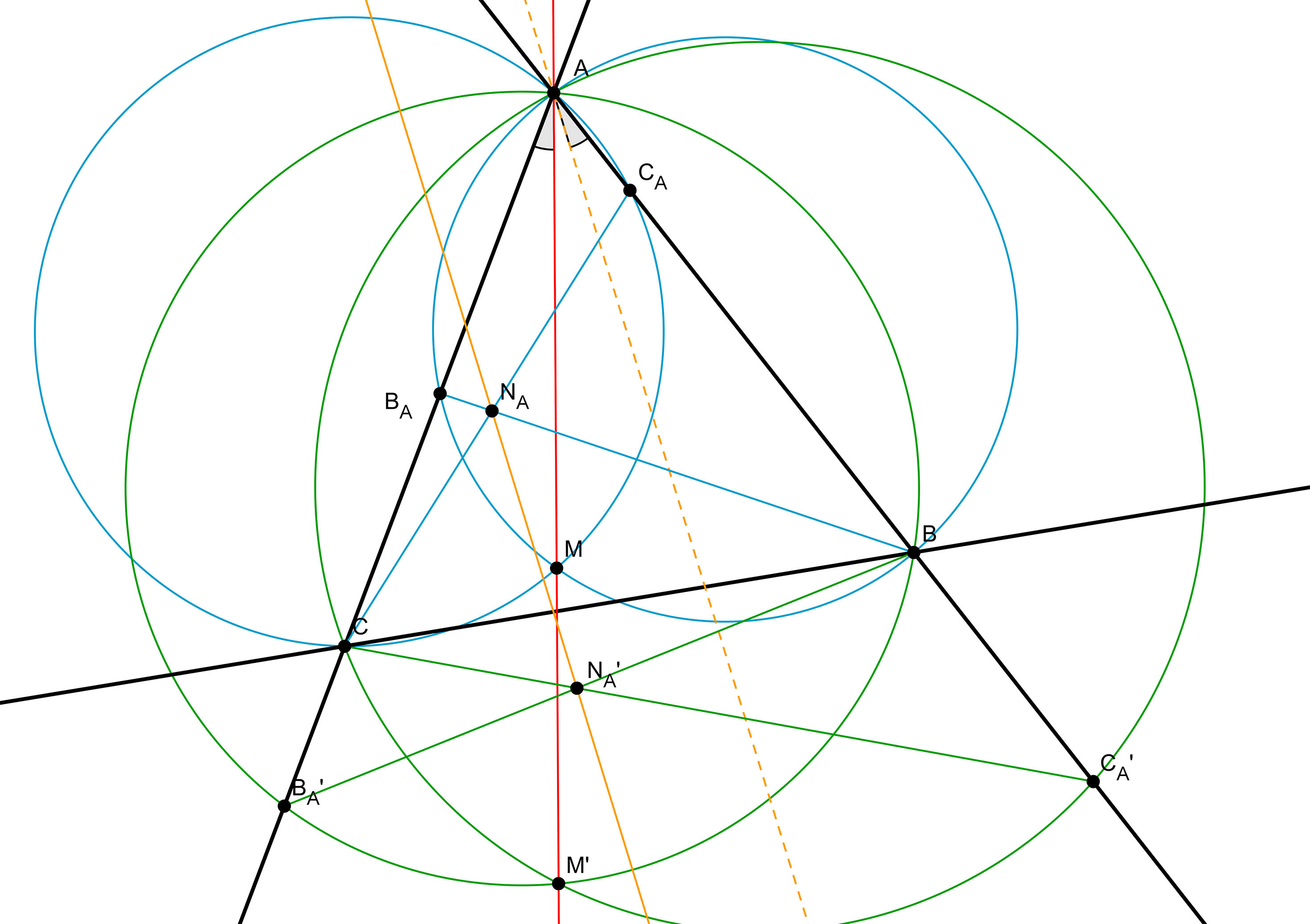}
			\newline
			Fig.19~Statement 4
		\end{center}
		\begin{proof}
			Strong proof can be provided with using methods of complex plane calculations. 
			Let us consider the following 
			\begin{eqnarray*}
				& A=0\in\mathbb{C}\\
				& B=b_0e^{i\beta_0}\in\mathbb{C}\\
				& C=1\in\mathbb{C}\\
				& M_A=me^{i\mu_0}, m\in\mathbb{R}\setminus\left\{{0}\right\}
			\end{eqnarray*}
			Here $m$ is free parameter of line $l$. $m\neq{0}$ because $M_A\neq{A}$.
			\newline
			We need to prove that 
			\begin{equation*}
				N_A=n_0+ne^{i\left(\beta_0-\mu_0\right)},~n\in\mathbb{R}
			\end{equation*}
			It means that $N_A$ belongs to line with parameter $n$ and this line is parallel to isogonal of $l$. 
			\newline
			For lines $AB, AC$ we have
			\begin{eqnarray*}
				& AB{\ni}z=te^{i\beta_0},~t\in\mathbb{R}\\
				& AC{\ni}z=t\in\mathbb{R}\\
			\end{eqnarray*}
			And for perpendicular bisectors of line segments $AM, AB, AC$
			\begin{eqnarray*}
				& p_{AB}{\ni}z=tie^{i\beta_0}+\frac{b_0}{2}e^{i\beta_0},~t\in\mathbb{R}\\
				& p_{AC}{\ni}z=ti+\frac{1}{2},~t\in\mathbb{R}\\
				& p_{AM}{\ni}z=tie^{i\mu_0}+\frac{m}{2}e^{i\mu_0},~t\in\mathbb{R}
			\end{eqnarray*}
			Thus if $O_{ABM}, O_{ACM}$ are centres of circumcircles $\omega_{ABM},\omega_{ACM}$ then we can calculate their coordinates. 
			\begin{equation*}
				O_{ABM}=z~\Leftrightarrow~
				\begin{cases}
					z=t_1ie^{i\beta_0}+\frac{b_0}{2}e^{i\beta_0}=t_2ie^{i\mu_0}+\frac{m}{2}e^{i\mu_0}\\
					\overline{z}=-t_1ie^{-i\beta_0}+\frac{b_0}{2}e^{-i\beta_0}=-t_2ie^{-i\mu_0}+\frac{m}{2}e^{-i\mu_0}
				\end{cases} 
			\end{equation*}
			Or
			\begin{equation*}
				\begin{cases}
					ie^{i\beta_0}t_1-ie^{i\mu_0}t_2=\frac{m}{2}e^{i\mu_0}-\frac{b_0}{2}e^{i\beta_0}\\
					-ie^{-i\beta_0}t_1+ie^{-i\mu_0}t_2=\frac{m}{2}e^{-i\mu_0}-\frac{b_0}{2}e^{-i\beta_0}
				\end{cases}
			\end{equation*}
			After solving system
			\begin{equation*}
				t_2=\frac{m\operatorname{cos}\left(\mu_0-\beta_0\right)-b_0}{2\operatorname{sin}\left(\mu_0-\beta_0\right)}
			\end{equation*}
			And 
			\begin{equation*}
				O_{AMB}=\frac{m\operatorname{cos}\left(\mu_0-\beta_0\right)-b_0}{2\operatorname{sin}\left(\mu_0-\beta_0\right)}ie^{i\mu_0}+\frac{m}{2}e^{i\mu_0}
			\end{equation*}
			Analogously
			\begin{equation*}
				O_{AMC}=\frac{m\operatorname{cos}\mu_0-1}{2\operatorname{sin}\mu_0}ie^{i\mu_0}+\frac{m}{2}e^{i\mu_0}
			\end{equation*}
			Let us calculate $B_A=\omega_{AMB}\cap{AC}$
			\begin{equation*}
				\begin{cases}
					B_A\in{AC}~\Leftrightarrow~B_A=t\\
					B_A\in\omega_{AMB}~\Leftrightarrow~B_A\overline{B_A}=B_A\overline{O_{AMB}}+\overline{B_A}O_{AMB}
				\end{cases}
			\end{equation*}
			$B_A\neq{A}$. It follows that $t\neq{0}$. Also $t\in\mathbb{R}$. Therefore
			\begin{eqnarray*}
				&t=O_{AMB}+\overline{O_{AMB}}=\\
				&=\frac{m\operatorname{cos}\left(\mu_0-\beta_0\right)-b_0}{\operatorname{sin}\left(\mu_0-\beta_0\right)}\left(\frac{i}{2}e^{i\mu_0}-\frac{i}{2}e^{-i\mu_0}\right)+\frac{m}{2}\left(e^{i\mu_0}-e^{-i\mu_0}\right)
			\end{eqnarray*}
			It follows that
			\begin{equation*}
				B_A=m\operatorname{cos}\mu_0-\frac{m\operatorname{cos}\left(\mu_0-\beta_0\right)-b_0}{\operatorname{sin}\left(\mu_0-\beta_0\right)}\operatorname{sin}\mu_0
			\end{equation*}
			Or 
			\begin{equation*}
				B_A=m\frac{\operatorname{cos}\mu_0\operatorname{sin}\left(\mu_0-\beta_0\right)-\operatorname{cos}\left(\mu_0-\beta_0\right)\operatorname{sin}\mu_0}{\operatorname{sin}\left(\mu_0-\beta_0\right)}+
				b_0\frac{\operatorname{sin}\mu_0}{\operatorname{sin}\left(\mu_0-\beta_0\right)}
			\end{equation*}
			Finally
			\begin{equation*}
				B_A=\frac{b_0\operatorname{sin}\mu_0-m\operatorname{sin}\beta_0}{\operatorname{sin}\left(\mu_0-\beta_0\right)}
			\end{equation*}
			Let us calculate $C_A=\omega_{AMC}\cap{AB}$
			\begin{equation*}
				\begin{cases}
					C_A\in{AB}~\Leftrightarrow~C_A=te^{i\beta_0}\\
					C_A\in\omega_{AMC}~\Leftrightarrow~C_A\overline{C_A}=C_A\overline{O_{AMC}}+\overline{C_A}O_{AMC}
				\end{cases}
			\end{equation*}
			$C_A\neq{A}$. It follows that $t\neq{0}$. Also $t\in\mathbb{R}$. Therefore
			\begin{equation*}
				t=O_{AMC}e^{-i\beta_0}+\overline{O_{AMC}}e^{i\beta_0}=
			\end{equation*}
			\begin{eqnarray*}
				&=\frac{m\operatorname{cos}\mu_0-1}{\operatorname{sin}\mu_0}ie^{i\mu_0}\left(\frac{i}{2}e^{i\left(\mu_0-\beta_0\right)}-
				\frac{i}{2}e^{-i\left(\mu_0-\beta_0\right)}\right)+\\
				&+\frac{m}{2}\left(e^{i\left(\mu_0-\beta_0\right)}+e^{-i\left(\mu_0-\beta_0\right)}\right)
			\end{eqnarray*}
			It follows that
			\begin{equation*}
				C_A=m\operatorname{cos}\left(\mu_0-\beta_0\right)e^{i\beta_0}-\frac{m\operatorname{cos}\mu_0-1}{\operatorname{sin}\mu_0}\operatorname{sin}\left(\mu_0-\beta_0\right)e^{i\beta_0}
			\end{equation*}
			Or
			\begin{equation*}
				C_A=\left(
				m\frac{\operatorname{sin}\mu_0\operatorname{cos}\left(\mu_0-\beta_0\right)-\operatorname{cos}\mu_0\operatorname{sin}\left(\mu_0-\beta_0\right)}{\operatorname{sin}\mu_0}+\frac{\operatorname{sin}\left(\mu_0-\beta_0\right)}{\operatorname{sin}\mu_0}
				\right)e^{i\beta_0}
			\end{equation*}
			Finally
			\begin{equation*}
				\label{eq:eq14}
				C_A=\frac{m\operatorname{sin}\beta_0+\operatorname{sin}\left(\mu_0-\beta_0\right)}{\operatorname{sin}\mu_0}e^{i\beta_0}
			\end{equation*}
			Let $p, q$ are such that
			\begin{eqnarray*}
				&B_A=p\\
				&C_A=qe^{i\beta_0}
			\end{eqnarray*}
			Therefore we have equations for line $BB_A$ and line $CC_A$ 
			\begin{equation*}
				BB_A{\ni}z=\left(b_0e^{i\beta_0}-p\right)t+b_0e^{i\beta_0}, t\in\mathbb{R}
			\end{equation*}
			\begin{equation*}
				CC_A{\ni}z=\left(1-qe^{i\beta_0}\right)t+1, t\in\mathbb{R}
			\end{equation*}
			It follows that calculating of point $N_A=BB_A\cap{}CC_A$ is solving following system
			\begin{equation*}
				\begin{cases}
					z=\left(1-qe^{i\beta_0}\right)t_1+1=\left(b_0e^{i\beta_0}-p\right)t_2+b_0e^{i\beta_0}\\
					\overline{z}=\left(1-qe^{-i\beta_0}\right)t_1+1=\left(b_0e^{-i\beta_0}-p\right)t_2+b_0e^{-i\beta_0}
				\end{cases}
			\end{equation*}
			Or 
			\begin{equation*}
				\begin{cases}
					\left(1-qe^{i\beta_0}\right)t_1-\left(b_0e^{i\beta_0}-p\right)t_2=b_0e^{i\beta_0}-1\\
					\left(1-qe^{-i\beta_0}\right)t_1-\left(b_0e^{-i\beta_0}-p\right)t_2=b_0e^{-i\beta_0}-1
				\end{cases}
			\end{equation*}
			After solving system 
			\begin{equation*}
				t_1=b_0\frac{p-1}{b_0-pq}
			\end{equation*}
			And 
			\begin{equation*}
				N_A=\left(1-qe^{i\beta_0}\right)t_1+1
			\end{equation*}
			For calculating $t_1$ we need to calculate $(p-1)$ and $b_0-pq$. Therefore
			\begin{equation*}
				p-1=\frac{b_0\operatorname{sin}\mu_0-m\operatorname{sin}\beta_0}{\operatorname{sin}\left(\mu_0-\beta_0\right)}-1=
				\frac{
				b_0\operatorname{sin}\mu_0-m\operatorname{sin}\beta_0-\operatorname{sin}\left(\mu_0-\beta_0\right)
				}{\operatorname{sin}\left(\mu_0-\beta_0\right)}
			\end{equation*}
			And 
			\begin{equation*}
				b_0-pq=b_0-\left(
				\frac{b_0\operatorname{sin}\mu_0-m\operatorname{sin}\beta_0}{\operatorname{sin}\left(\mu_0-\beta_0\right)}
				\right)
				\left(
				\frac{m\operatorname{sin}\beta_0+\operatorname{sin}\left(\mu_0-\beta_0\right)}{\operatorname{sin}\mu_0}
				\right)
				=
			\end{equation*}
			\begin{equation*}
				=-m\operatorname{sin}\beta_0
				\frac{b_0\operatorname{sin}\mu_0-m\operatorname{sin}\beta_0-\operatorname{sin}\left(\mu_0-\beta_0\right)}{\operatorname{sin}\mu_0\operatorname{sin}\left(\mu_0-\beta_0\right)}
			\end{equation*}
			It follows that
			\begin{equation*}
				b_0-pq=-\frac{(p-1)m\operatorname{sin}\beta_0}{\operatorname{sin}\mu_0}
			\end{equation*}
			Thus 
			\begin{equation*}
				t_1=-b_0\frac{\operatorname{sin}\mu_0}{m\operatorname{sin}\beta_0}
			\end{equation*}
			Also
			\begin{equation}
				\label{eq:eq13}
				N_A=\left(1-qe^{i\beta_0}\right)t_1+1=1+t_1-qt_1\operatorname{cos}\beta_0-iqt_1\operatorname{sin}\beta_0
			\end{equation}
			Let us noticed that 
			\begin{equation*}
				\operatorname{sin}\mu_0-\operatorname{sin}\left(\mu_0-\beta_0\right)\operatorname{cos}\beta_0=
				\operatorname{sin}\beta_0\operatorname{cos}\left(\mu_0-\beta_0\right)
			\end{equation*}
			It follows that 
			\begin{eqnarray*}
				&t_1(1-q\operatorname{cos}\beta_0)=\\
				&=-b_0
				\left(
				\frac{\operatorname{sin}\mu_0}{m\operatorname{sin}\beta_0}
				\right)
				\left(
				\frac{\operatorname{sin}\mu_0-\operatorname{sin}\left(\mu_0-\beta_0\right)\operatorname{cos}\beta_0
				-m\operatorname{sin}\beta_0\operatorname{cos}\beta_0}{\operatorname{sin}\mu_0}
				\right)
				=
			\end{eqnarray*}
			\begin{equation*}
				=b_0\operatorname{cos}\beta_0-\frac{b_0}{m}\operatorname{cos}\left(\mu_0-\beta_0\right)
			\end{equation*}
			With (\ref{eq:eq13}) we have
			\begin{equation}
				\label{eq:eq14}
				\operatorname{Re}N_A=1+b_0\operatorname{cos}\beta_0-\frac{b_0}{m}\operatorname{cos}\left(\mu_0-\beta_0\right)
			\end{equation}
			Also
			\begin{eqnarray*}
				&qt_1\operatorname{sin}\beta_0=-b_0
				\left(
				\frac{\operatorname{sin}\mu_0}{m\operatorname{sin}\beta_0}
				\right)
				\left(
				\frac{m\operatorname{sin}\beta_0+\operatorname{sin}\left(\mu_0-\beta_0\right)}{\operatorname{sin}\mu_0}
				\operatorname{sin}\beta_0
				\right)=\\
				&=-b_0\operatorname{sin}\beta_0-\frac{b_0}{m}\operatorname{sin}\left(\mu_0-\beta_0\right)
			\end{eqnarray*}
			With (\ref{eq:eq13}) we have
			\begin{equation}
				\label{eq:eq15}
				\operatorname{Im}N_A=b_0\operatorname{sin}\beta_0+\frac{b_0}{m}\operatorname{sin}\left(\mu_0-\beta_0\right)
			\end{equation}
			Finally with (\ref{eq:eq13})--(\ref{eq:eq15}) and $n=\frac{1}{m},~n\in\mathbb{R}$ we have
			\begin{equation*}
				N_A=1+e^{i\beta_0}-nb_0e^{i\left(\beta_0-\mu_0\right)}
			\end{equation*}
		\end{proof}
		We can give also converse proposition.
		\begin{corollary}
			\label{cor:cor2}
			If cevians intersection point $N_A\neq{A}$ belongs to fixed line $l$ and $A\in{l}$ then respective Miquel-Steiner's point $M_A$ belongs to circle which is isogonal of $l$. (Fig.20)
		\end{corollary}
		\begin{center}
			\includegraphics[width=0.5\textwidth]{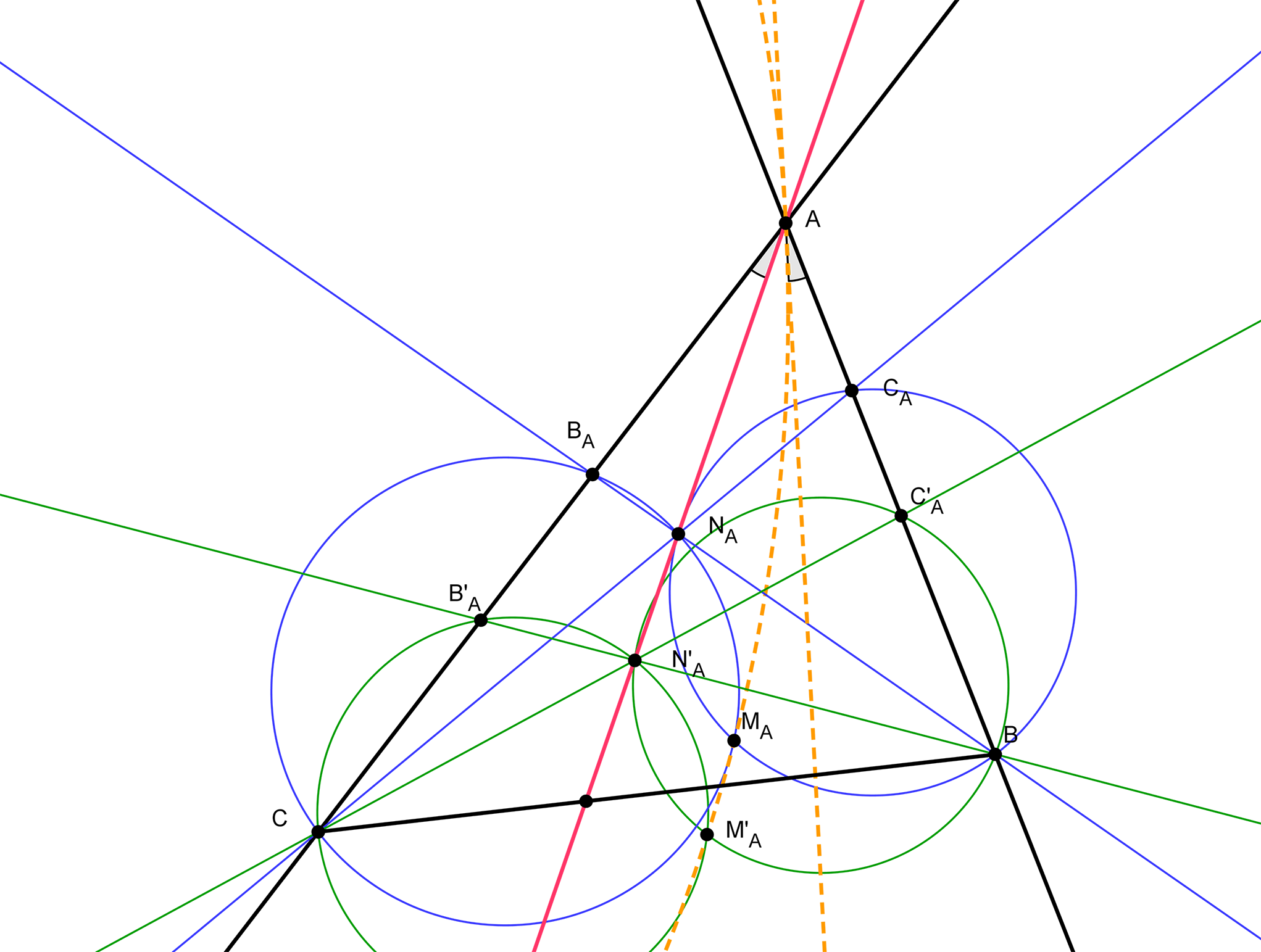}
			\newline
			Fig.20~Corollary 2
		\end{center}
		\begin{proof}
			Proof is simple. We just need to fix some circle with centre $A$ and use inversion with respect to this circle [5, p.78]. With this inversion we will have construction of statement \ref{st:st4}. We have needed because of inversion is involutive and conformal.
		\end{proof}
	\section*{Miquel-Steiner's point and triangle centres}
		Here we consider following situations. Miquel-Steiner's point is 
		\newline 1)~incenter, 
		\newline 2)~orthocenter, 
		\newline 3)~circumcenter.
		\begin{statement}
			\label{st:st5}
			If Miquel-Steiner's point $M_A=I$ is incenter then respective cevians $BB_A, CC_A$ are such that $\left|BC_A\right|=\left|CB_A\right|=\left|BC\right|$. (Fig.21)
		\end{statement}
		\begin{center}
			\includegraphics[width=0.5\textwidth]{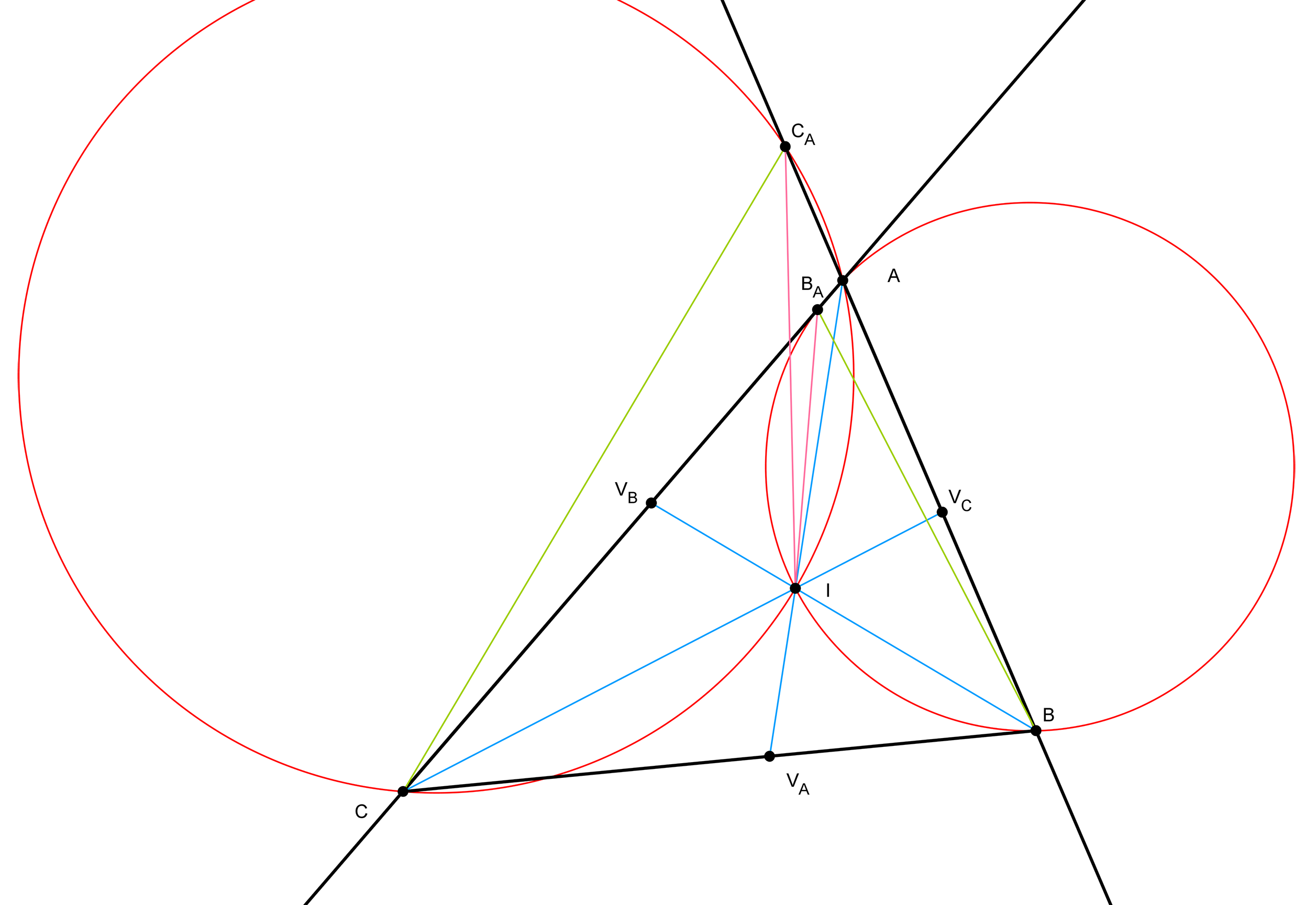}
			\newline
			Fig.21~Miquel-Steiner's point is incenter
		\end{center}
		\begin{proof}
			Let $V_A=AI\cap{BC}$ -- bisector base. Analogously $V_B, V_C$. It is well-known that equal angles gives equal chordes. Thus $\left|B_AI\right|=\left|BI\right|$. Therefore we have angles equalities.
			\begin{eqnarray*}
				& \angle{IBB_A}=\angle{BB_AI} \\
				& \angle{IB_AC}=\angle{V_BBA}=\angle{CBV_B}\\
			\end{eqnarray*}
			Thus $\angle{BB_AC}=\angle{CBB_A}$. Or $\left|CB_A\right|=\left|BC\right|$. Analogously $\left|BC_A\right|=\left|BC\right|$.
			\qedhere
		\end{proof}
		\begin{remark}
			\label{rm:rm5}
			Converse proposition gives us four probable Miquel-Steiner's points. They are intesection points of bisectors (internal and external).
		\end{remark}
		\begin{statement}
			\label{st:st6}
			Miquel-Steiner's point $M_A=H$ is orthocenter if and only if respective cevians $BB_A, CC_A$ are such that $\left|BB_A\right|=\left|CC_A\right|=\left|BC\right|$. (Fig.22)
		\end{statement}
		\begin{center}
			\includegraphics[width=0.5\textwidth]{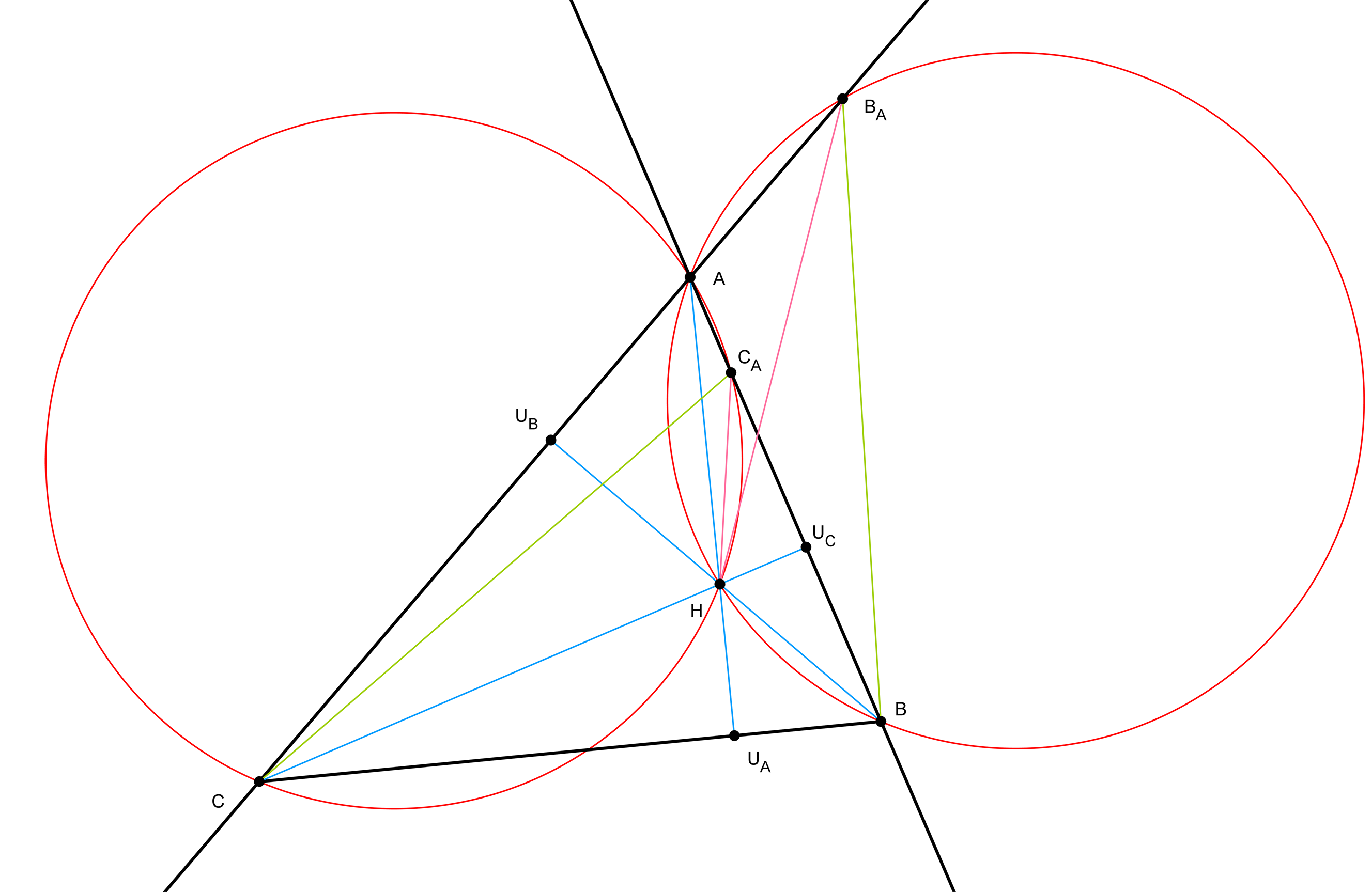}
			\newline
			Fig.22~Miquel-Steiner's point is orthocenter
		\end{center}
		\begin{proof}
			Let us proof if Miquel-Steiner's point $M_A=H$ is orthocenter then  $\left|BB_A\right|=\left|CC_A\right|=\left|BC\right|$. Converse proposition is provided because of unique pair of cevians such that $\left|BB_A\right|=\left|CC_A\right|=\left|BC\right|$.
			\newline
			Let $U_A=AH\cap{BC}$ -- altitude base. Analogously $U_B, U_C$. It is obvious that
			\begin{equation*}
				\angle{U_BBB_A}=\angle{HAC}=90^{\circ}-\angle{ACB}=\angle{CBU_B}
			\end{equation*}
			Thus $BU_B$ is bisector and altitude of triangle $\Delta{BB_AC}$. It follows it is isosceles and $\left|BB_A\right|=\left|BC\right|$. Analogously $\left|CC_A\right|=\left|BC\right|$.
			\qedhere
		\end{proof}
		\begin{statement}
			\label{st:st7}
			Miquel-Steiner's point $M_A=O$ is circumcenter if and only if respective cevians $BB_A, CC_A$ are such that $\left|BB_A\right|=\left|CB_A\right|, \left|CC_A\right|=\left|BC_A\right|$. (Fig.23) Moreover $N_A\in{\omega_{ABC}}$.
		\end{statement}
		\begin{center}
			\includegraphics[width=0.5\textwidth]{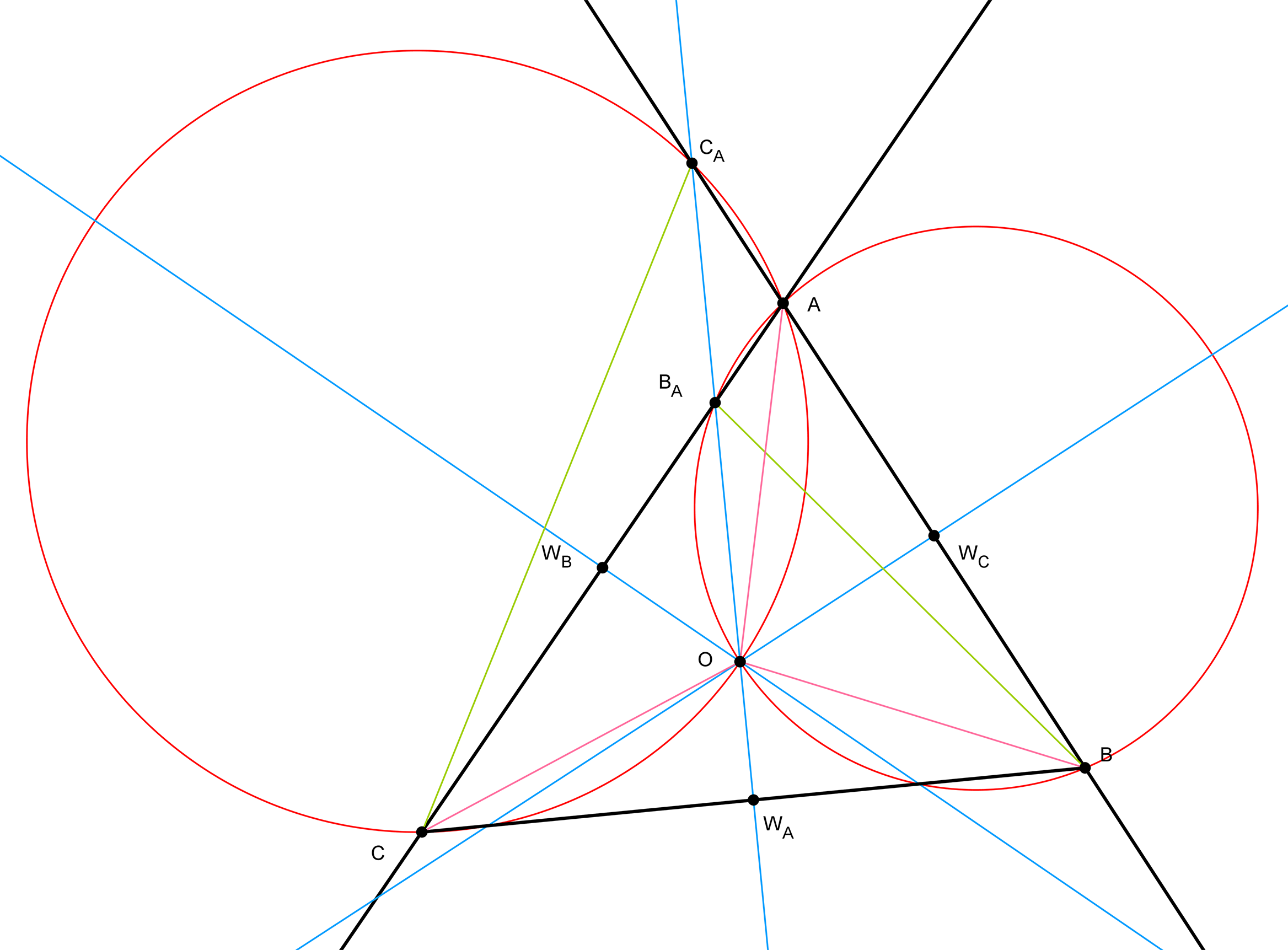}
			\newline
			Fig.23~Miquel-Steiner's point is circumcenter
		\end{center}
		\begin{proof}
			Let us proof if Miquel-Steiner's point $M_A=O$ is circumcenter then $\left|BB_A\right|=\left|CB_A\right|$, $\left|CC_A\right|=\left|BC_A\right|$. Converse proposition is provided because of unique pair of cevians such that $\left|BB_A\right|=\left|CB_A\right|$, $\left|CC_A\right|=\left|BC_A\right|$.
			\newline
			It is obvious that
			\begin{equation*}
				\angle{BB_AO}=\angle{BAO}=\angle{OBA}=\angle{OB_AC}
			\end{equation*}
			Also 
			\begin{equation*}
				\angle{ACO}=\angle{OAC}=\angle{OBB_A}
			\end{equation*}
			But $\angle{OCB}=\angle{COB}$. Therefore $\left|BB_A\right|=\left|CB_A\right|$. Analogously $\left|CC_A\right|=\left|BC_A\right|$. 
			\newline
			Now we have two isosceles triangles $\Delta{BB_AC}, \Delta{BC_AC}$ with similar base $BC$. Thus they have colinear bisectors. It follows that $B_A, C_A, O$ are colinear. Therefore we have
			\begin{equation*}
				\angle{B_ABA}=\angle{CBA}-\angle{CBB_A}=\angle{C_ACB}-\angle{B_ACB}=\angle{C_ACA}
			\end{equation*}
			It follows that $BB_A, CC_A$ intersects on $\omega_{ABC}$.
			\qedhere
		\end{proof}
	\section*{Conclusion}
		We have seen that the Miquel-Steiner's theorem can be reformulated for triangle $\Delta{ABC}$ and two cevians $BB_A, CC_A$. With this point of view we can stydy Miquel-Steiner's point locus. Firstly case of parallel cevians was studied (lemma \ref{lm:lm1}). Then for case when cevians intersection point exists and unique the bijection(or one-to-one correspondence) between pair of cevians and the Miquel-Steiner's point was constructed (statement \ref{st:st1} and corollary \ref{cor:cor1}). In this case Miquel-Steiner's point locus is entire plane without lines $AB, BC$ and circumcircle $\omega_{ABC}$. One of the main result of this article is theorem \ref{th:th1}. It is about Miquel-Steiner's point locus for internal cevians. With this theorem for fixed vertex $A$ Miquel-Steiner's auxiliary circles, auxiliary and main centres, axis were defined (definition \ref{def:def1}). For these construction few properties were obtained and proved in 2012 as a part of Miquel-Steiner's point study (theorems \ref{th:th2}--\ref{th:th4} and lemma \ref{lm:lm2}). Then it was found that these results had be obtained earlier by Henri Brocard (remarks \ref{rm:rm2}--\ref{rm:rm3}). 
		\newline
		Also few results of specific cases were obtained and proved. Statement \ref{st:st2} is related to perpendicular cevians, statement \ref{st:st3} is related to cevians are such that $B_AC_A \mid \mid BC$. Remark \ref{rm:rm4} gives us understanding that cevians intersection point lying on median corresponds to the Miquel-Steiner's point lying on symmedian. Statement \ref{st:st4} is logical extension of remark \ref{rm:rm4}. It is related to situation when the Miquel-Steiner's point belong to some line from $A$. Also we gave converse proposition (corollary \ref{cor:cor2}).
		\newline
		Also we studied cases when Miquel-Steiner's point is coincides with incenter (statement \ref{st:st5} and remark \ref{rm:rm5}), orthocenter (statement \ref{st:st6}) and circumcenter (statement \ref{st:st7}).

\end{document}